\documentclass{article}
\usepackage{mathrsfs,amsmath, amssymb}
\usepackage{amsthm}
\usepackage{graphicx}
\usepackage{caption}
\usepackage[labelformat=simple]{subcaption}

\DeclareCaptionLabelFormat{carinasub}{#2.}
\theoremstyle{remark}
\newtheorem*{remarks}{Remarks}
\usepackage{nicefrac}
\usepackage{tikz}
\usepackage{authblk}
\usepackage{bbm}
\usepackage{float}
\usepackage{stmaryrd}
\usepackage{biblatex} 
\usepackage{hyperref}
\usepackage{cleveref}
\usepackage{nicefrac,xfrac}
\usepackage{algorithm}
\usepackage{algpseudocode}
\usepackage{xcolor}
\usepackage{subcaption}
\usepackage{caption}
\usepackage[labelformat=simple]{subcaption}

\DeclareCaptionLabelFormat{carinasub}{#2.}
\makeatletter
\renewcommand{\p@subfigure}{}
\makeatother
\numberwithin{equation}{section}
\providecommand{\tol}[1]{\operatorname{tol}_{\text{#1}}}

\theoremstyle{remark}

\title{ \bf{Adaptive FEM with explicit time integration for the wave equation}}
\author[1]{Marcus J. Grote}
\author[2]{Omar Lakkis}
\author[1]{Carina S. Santos}
\affil[1]{Department of Mathematics and Computer Science, University of Basel, Spiegelgasse 1,
4051 Basel, Switzerland}
\affil[2]{Department of Mathematics, University of Sussex, 
Brighton, England, UK}
\date{\today}

\begin{document}

\maketitle

\begin{abstract}
Starting from a recent a posteriori error estimator for
the finite element solution of the wave equation with explicit
time-stepping [Grote, Lakkis, Santos, 2024], we devise a space-time adaptive strategy which includes 
both time evolving meshes and local time-stepping [Diaz, Grote, 2009] to overcome any overly
  stringent CFL stability restriction on the time-step due to local mesh
  refinement. Moreover, at each time-step the adaptive algorithm monitors the 
  accuracy thanks to the error indicators and recomputes the current 
  step on a refined mesh until the desired tolerance is met; meanwhile, the mesh is coarsened
  in regions of smaller errors. Leapfrog based local time-stepping is applied in all regions of local mesh refinement to incorporate adaptivity into fully explicit time
  integration with mesh change while retaining efficiency.
  Numerical results illustrate the optimal rate of convergence 
  of the a posteriori error estimators on time evolving meshes.
\end{abstract}

\section{Introduction}
Numerical methods for the solution of PDEs aim to compute with utmost efficiency an approximation $u_h$ of the true solution $u$ for a prescribed error tolerance. Finite element methods (FEMs), be they continuous or discontinuous, offer a wide range of applicability by accommodating spatially varying media in complex geometry and discontinuous material interfaces. When combined with an adaptive mesh strategy, which concentrates the degrees of freedom required to accurately represent $u_h$ only where needed while keeping the computational cost small elsewhere, adaptive FEMs improve the computational efficiency and accuracy of numerical methods for wave propagation. This involves dynamically refining or moving the mesh in regions where the solution changes rapidly, or where errors are most significant, leading to improved accuracy with fewer degrees of freedom. 

A posteriori error estimates are the cornerstone of any adaptive strategy with guaranteed rigorous error bounds.
For elliptic problems, a posteriori error analysis is well-developed and
leads to rigorous and explicitly computable error bounds, which allow the end user to assess the accuracy of $u_h$, see \cite{AinsworthOden:00:book:A-posteriori,Verfurth:13:book:A-posteriori} and references therein. Typically, such error bounds consist of local contributions, known as error indicators, which can be used to automatically steer local mesh refinement strategy and thus improve upon the accuracy of $u_h$ while keeping the added cost minimal. For parabolic problems, both time-discretization and mesh change must also be included in the a posteriori error estimates. Various a posteriori error bounds are available for parabolic problems either based on duality or space-time Galerkin formulation, often coupled with a discontinuous Galerkin (DG) formulation in time \cite{ErikssonJohnson:91:article:Adaptive,
  Picasso:98:article:Adaptive,Verfurth:03:article:A-posteriori, 
  ChenJia:04:article:An-adaptive, Verfurth:13:book:A-posteriori,
  GaspozSiebertKreuzerZiegler:19:article:A-convergent}.

For (time-dependent) wave equations, error estimation and adaptivity are far less developed than for elliptic or parabolic problems. 
Indeed, in contrast to parabolic problems, the most commonly used time integration methods for wave equations, such as the popular second-order leapfrog (LF) method (or St\"ormer-Verlet), are explicit. Explicit time integration is highly efficient, in particular for large-scale problems on massively parallel architectures. In the presence of local mesh refinement, however, the CFL stability constraint will impose a tiny time-step across the entire computational domain which will cripple any explicit time integrator; hence,
standard explicit methods are generally deemed inefficient 
when combined with adaptivity as
''for adaptive methods to time-dependent
waves, unconditionally stable time-stepping methods are needed'' (\cite{ThompsonHe:05:article}, p. 1948).

It is no wonder thus that most a posteriori error estimates and space-time adaptive strategies for the wave equation have relied on implicit time integration
\cite{Johnson:93:article:Discontinuous, HulbertHughes:90:article:Space-Time-Finite,BernardiSuli:05:article:Time}.
Residual based a posteriori error estimates with first-order implicit time-stepping
were developed in \cite{BernardiSuli:05:article:Time,Adjerid:02:article:A-posteriori,
  Adjerid:06:article:A-posteriori}.
Alternatively, goal-oriented adaptivity always requires the solution of an adjoint
(dual) problem \cite{BangerthRannacher:01:article:Adaptive, BangerthGroteHohenegger:04:article}. Recently, various 
a posteriori error estimates were derived for semi-discrete
formulations with anisotropic mesh refinement either in the
$L^2(0,T; H^1 (\Omega)$-norm
\cite{Picasso:10:article:Numerical,%
  GoryninaLozinskiPicasso:19:article:An-Easily-Computable}
or in a ``damped energy norm''  \cite{Chaumont-Frelet:23:article:Asymptotically-Constant-Free}; the latter recently led to fully discrete error estimates
 \cite{Chaumont-FreletErn:24:techreport:Damped-Energy-Norm} yet under 
 a fixed mesh assumption. Unconditionally stable space-time Galerkin discretization of the second-order wave equation were also derived in \cite{ferrari2025intrinsicunconditionalstabilityspacetime}
 and, for star-shaped spatial domains, in \cite{BignardiMoiola:25:article:A-Space--Time-Continuous}
 .
 In \cite{ZankSteinbach:16:article:Adaptive-Space-Time}, 
 an adaptive algorithm was proposed for a boundary-element based space-time discretization in one space dimension without mesh-change.
Recently, a space-time adaptive algorithm was proposed for the wave equation using an implicit time-stepping strategy \cite{DongGeorgoulisMascottoWang:25:techreport:A-Posteriori-Error}.
 
In the design of space-time adaptive methods, a posteriori error estimates of fully discrete formulations in time-dependent problems need to include the added effect on the error due to mesh change 
from one time-step to the next;
in fact, some of the above cited
works address that particular issue for parabolic problems -- see also \cite{Dupont:82:article:Mesh,LakkisPryer:12:article:Gradient,SchmidtSiebert:05:book:Design}. 
For both accuracy and efficiency, it is indeed expected, often even required, 
for an adaptive method to locally refine and coarsen
the mesh repeatedly during the entire simulation.
Although quantification of mesh-change error in second order
hyperbolic problems is less studied, a notable result in this
direction was provided by
\cite{KarakashianMakridakis:05:article:Convergence} in an a priori
setting. 



To overcome the above mentioned difficulties, we recently derived a posteriori
error estimates for the wave equation with explicit LF time integration
and mesh change \cite{GroteLakkisSantos:24:PosterioriErrorEstimatesWave}, 
starting from the a posteriori analysis for
semi-discrete formulations \cite{GeorgoulisLakkisMakridakisVirtanen:16:article:A-Posteriori} (continous in space) 
discretized in time using a two-step Newmark (or cosine) family of 
methods, which include the explicit leapfrog method.
Moreover, our estimates also allow for LF based explicit local time-stepping (LF-LTS), which circumvents the local CFL bottleneck by taking smaller local time-steps, but only where needed, without sacrificing explicitness. 

Here we propose an adaptive strategy for the solution of the time-dependent wave equation based on rigorous a posteriori error estimates which include both mesh change
and LF based fully explicit local time-stepping. In Section 2, we present the finite element discretization including the Galerkin FE formulation with explicit LF time integration and mesh change.
In Section 3, 
starting from the a posteriori error estimators from \cite{GroteLakkisSantos:24:PosterioriErrorEstimatesWave}, which we briefly recall in \ref{the:full-error analysis}, we develop a new adaptive algorithm for the numerical solution
of the wave equation which combines three essential ingredients: a posteriori estimates for guaranteed error bounds, mesh change to adapt the mesh from one time-step to another, and local time-stepping to retain the efficiency of explicit time integration even in the presence of local mesh refinement.
In fact, local error indicators not only permit one to estimate but also to control the effect of mesh change and adaptivity and thus to decide ``on the fly'' whether to accept or recompute the numerical solution at the current time-step. By restricting
mesh change to compatible meshes, as detailed in Section 3.2, we ensure that local refinement never increases the numerical error while keeping under control the inherent but small error from mesh coarsening. Finally, in Section 4, we present 
four numerical experiments in one and two space dimensions
, which
confirm the expected convergence order and illustrate the usefulness of
our adaptive strategy.

\section{Galerkin FE formulation with explicit time integration and mesh change}
We consider the wave equation in a bounded Lipschitz domain $\Omega\subset \mathbbmss{R}^{d}$ 
\begin{equation}
\label{eq:waveeq}
\frac{\partial^2 }{\partial t^2}u(x, t)-\nabla \cdot\left[c^2(x)\nabla u(x, t)\right]=f(x, t), \qquad x \in \Omega, t \in(0, T],
\end{equation}
for a given wave speed $0<c_{\min} \leq c(x) \leq c_{\max}$, forcing $f(x,t)$ 
and initial conditions

\begin{equation}
u(x,0)=u_{0}(x), \qquad  \frac{\partial}{\partial t} u(x,0)=v_{0}(x), \qquad x \in \Omega.
\end{equation}
In addition, we impose for simplicity homogeneous Dirichlet or Neumann
boundary conditions at the boundary $\partial\Omega$ of $\Omega$
\begin{equation}
u(x,t)=0, \quad x\in\Gamma_0,\qquad \frac{\partial}{\partial n}u(x, t) = 0, \quad x\in\partial\Omega \setminus \Gamma_0, \qquad t > 0,
\end{equation}
where the Dirichlet boundary, $\Gamma_0$, has positive measure, $|\Gamma_0|>0$.

For $u_0 \in H^1_{\Gamma_0}(\Omega)$, $v_0 \in L^2(\Omega)$ and $f\in L^2(0,T; L^2(\Omega))$, it is well-known that
the wave equation \eqref{eq:waveeq} in fact has a unique solution $u\in C^0(0,T;  H^1_{\Gamma_0}(\Omega)) \cap C^1(0,T;  L^2(\Omega))$ \cite{LionsMagenes:72:.Non-HomogeneousBoundaryValueProblemsandApplications}.

First, we discretize time with a fixed time-step $\Delta t>0$ on a standard uniform
time-grid with integer indices
\begin{equation}
0=t_{0}<t_{1}<\cdots<t_{N}=T, \qquad t_{n}=n \Delta t .
\end{equation}

%
%

To each time $t_{n}, n=0, \ldots, N$, we then associate a spatial mesh $\mathcal{M}_{n}$ 
and a corresponding finite element space
\begin{equation}
\mathbbmss{V}_{n}=\mathbbmss{P}^{k}\left(\mathcal{M}_{n}\right) \cap  C^{0}(\Omega)
\end{equation}
of continuous piecewise polynomials of degree $k\geq 1$. The finite element mesh 
$\mathcal{M}_{n}$ consists of polytopal (triangular, quadrilateral, tetrahedral, etc.) elements $K \in \mathcal{M}_{n}$ with flat sides (edges, faces, etc.).
The corresponding piecewise constant \emph{meshsize} function
\begin{equation}
h^{n}(x):=\operatorname{diam} \bigcap_{\boldsymbol{x} \in K \in \mathcal{M}_{n}} \bar{K}.
\end{equation}
For each $E \in \mathcal{M}_{n} \cup \operatorname{Sides} \mathcal{M}_{n}$ we denote the constant value of $h^n$ by $h_{E}$ (which thus depends implicitly on $n$).

With adaptvity in mind, we allow for the case where $\mathcal{M}_{n}$ (and thus $\mathbbmss{V}_{n}$) changes with time under \emph{mesh compatibility conditions}, which implies that at each point of the domain either $\mathcal{M}_{n-1}$ is a 
(compatible) refinement of $\mathcal{M}_{n}$ or conversely -- see Section \ref{sec:compatible_mesh_change} and \cite{LakkisMakridakis:06:article:Elliptic} for details.

To each $\mathbbmss{V}_{n}$ we associate a Lagrangian nodal finite element basis 
\begin{equation}
\left\{\Phi_{1}^{n}, \ldots, \Phi_{M_{n}}^{n}\right\}, \qquad \Phi_{m}^{n}(z_\ell^n) = \delta_{m,\ell} \qquad M_{n}:=\operatorname{dim} \mathbbmss{V}_{n}
\end{equation}
and corresponding finite element nodes $z_{m}^{n}$, 
$m=1, \ldots, M_{n}$. 

To transfer functions between FE spaces during mesh change, for instance,
we introduce the \textit{FE transfer operator} $\Pi_{n}:  C^{0}(\Omega) \rightarrow \mathbbmss{V}_{n}$, which may either 
denote the $L^2$-projection or standard Lagrangian interpolation

\begin{equation}
\label{eq:Piop}
\Pi_{n} v(x)=\sum_{m=1}^{M_{n}} v(z_m^n)\Phi_{m}^{n}(x).
\end{equation}

%

Each mesh $\mathcal{M}_{n}$ has two types of elements, \emph{coarse} and \emph{fine}, 
$\mathcal{M}_{n}=\mathcal{M}_{n}^{\mathrm{c}} \cup \mathcal{M}_{n}^{\mathrm{f}}$, 
%
%
while a basis function $\Phi_{m}^{n}$ is \emph{fine} only if its support intersects at least one element in the fine mesh $\mathcal{M}_{n}^{\mathrm{f}}$; otherwise, it is \emph{coarse}. Letting $\mathbbmss{V}_{n}^{\mathrm{f}}$ and $\mathbbmss{V}_{n}^{\mathrm{c}}$ respectively be the subspaces spanned by the fine and coarse basis functions, we thus have $\mathbbmss{V}_{n}=\mathbbmss{V}_{n}^{\mathrm{f}} \oplus \mathbbmss{V}_{n}^{\mathrm{c}}$. If the indices are ordered into fine-first from $1, \ldots, M_{n}^{\mathrm{f}}$, for some integer $M_{n}^{\mathrm{f}} \leq M_{n}$, and coarse-last $M_{n}^{\mathrm{f}}+1, \ldots, M_{n}$ every finite element function $V\in\mathbbmss{V}_{n}$ can be written as
\begin{equation}
V(x)=\left(\sum_{m=1}^{M_{n}^{\mathrm{f}}}+\sum_{m=M_{n}^{\mathrm{f}}+1}^{M_{n}}\right) \Phi_{m}^{n}(x) \mathrm{v}_{m}
\end{equation}
for a suitable vector $\mathbf{v}=\left(\mathrm{v}_{1}, \ldots, \mathrm{v}_{M_{n}}\right) \in \mathbbmss{R}^{M_{n}}$. 
Similarly to \eqref{eq:Piop}, we define the \emph{fine-mesh interpolator}
$\Pi_{n}^{\mathrm{f}}:  C^{0}(\Omega) \rightarrow \mathbbmss{V}_{n}^{\mathrm{f}}$
as
\begin{equation}
\Pi_{n}^{\mathrm{f}} V=\sum_{m=1}^{M_{n}^{\mathrm{f}}} \Phi_{m}^{n} \mathrm{v}_{m}, \qquad V \in \mathbbmss{V}_{n}.
\end{equation}

%
%
%

Given the continuous elliptic operator $\mathcal{A}$ and corresponding bilinear form $a$
\begin{equation}
\langle\mathcal{A} u \mid v\rangle =a(u,v) = \int_{\Omega} c^2(x) \nabla u \cdot \nabla v\, dx, \qquad  \forall u, v \in H^1(\Omega)
\end{equation}
we introduce
for each $n$ the corresponding \emph{discrete elliptic operator} $A_{n} = A_{\mathbbmss{V}_{n}}$ on $\mathbbmss{V}_{n}$, and the source approximation
\begin{equation}
F^n = \Pi_n F(., t_n).
\end{equation}

Here for any conforming finite element
subspace $\mathbbmss{V} \subseteq H^1_{\Gamma_0}$, $A_{\mathbbmss{V}}$ is defined as
\begin{equation}
\begin{aligned}
A_{\mathbbmss{V}}: \mathbbmss{V} & \rightarrow \mathbbmss{V} \\
 \phi & \mapsto A_{\mathbbmss{V}} \phi:\left\langle A_{\mathbbmss{V}} \phi, v\right\rangle=\langle\mathcal{A} \phi \mid v\rangle = a(\phi, v) \quad \forall \; v \in \mathbbmss{V}.
\end{aligned}
\end{equation}

In fact, we can naturally extend $A_{\mathbbmss{V}}$ to any larger (finite or infinite dimensional) conforming subspace 
$\mathbbmss{V} \subseteq \mathbbmss{W} \subseteq H^1_{\Gamma_0}$ as 
\begin{equation}
\begin{array}{rlll}
A_{\mathbbmss{V}}:   \mathbbmss{W} & \rightarrow & \mathbbmss{V} \\
w & \mapsto & A_{\mathbbmss{V}} w
\end{array}
\end{equation}
thanks to Riesz's representation theorem via
\begin{equation}
\left\langle A_{\mathbbmss{V}} w, v\right\rangle=\langle\mathcal{A} v \mid v\rangle  = a(\phi, v) \quad \forall \; v \in \mathbbmss{V}.
\end{equation}
Alternatively we can think of $A_{\mathbbmss{V}}=\Pi_{\mathbbmss{V}} \mathcal{A}$, where $\Pi_{\mathbbmss{V}}$ denotes the $\mathrm{L}_{2}$ projection onto $\mathbbmss{V}$.

Whenever local time-stepping is required, we replace $A_n$ formally by $\widetilde{A_{n}}$,
for instance  -- see Section 3.3 for further details:
\begin{equation}
\label{eq:LTS_A}
\widetilde{A_{n}}:=A_{n}-\frac{\Delta t^{2}}{16} A_{n} \Pi_{n}^{\mathrm{f}} A_{n}.
\end{equation}

The particular instance of $\widetilde{A_{n}}$ in \eqref{eq:LTS_A} corresponds to the simplest situation with two local time-steps of size $\Delta t / 2$ each for each global time-step of size $\Delta t$. By letting $\widetilde{A_{n}}$ denote a generic perturbed bilinear form induced by local time-stepping, our analysis inherently encompasses situations with different coarse-to-fine time-step ratios, too, which may even change from one locally refined subregion to another across a single mesh. In fact, it even includes a hierarchy of locally refined regions, each associated with its own local time-step \cite{DiazGrote:15:article:Multi-Level-Explicit}

For time discretization, we use a standard centered second-order finite difference 
approximation. In doing so, we always transfer $U^n$ and $U^{n-1}$ first to
$\mathbbmss{V}_{n+1}$ to ensure that $U^{n-1}$, $U^n$ and $U^{n+1}$ all three belong to the same FE space. This yields the time-varying finite element space leapfrog scheme:
\begin{equation}
\label{eq:leapfrog}
\begin{aligned}
U^{0} & :=\Pi_0 u_{0} \\
U^{1} & :=\Pi_{1}\left[U^{0}+ \Pi_{0} v_{0}\, \Delta t+\left(F^{0}-\widetilde{A_{0}} U^{0}\right) \Delta t^2 \right] \\
U^{n+1} & :=\Pi_{n+1}\left[2\, U^{n}-\Pi_{n} U^{n-1}+\left(F^{n}-\widetilde{A_{n}} U^{n}\right) \Delta t^{2}\right], \: n=1, \ldots, N.
\end{aligned}
\end{equation}
If needed, the velocity $V^{n+1 / 2}\in\mathbbmss{V}_{n+1}$ at intermediate times can easily be obtained through
\begin{equation}
V^{n+1 / 2} := \frac{U^{n+1} - \Pi_{n+1} U^{n}}{\Delta t}.
\end{equation}

\section{Space-time adaptive finite element method}\label{sec:space-time}
We shall now develop a fully adaptive method for the numerical solution
of the wave equation with combines three essential ingredients: a posteriori estimates for guaranteed error bounds, mesh change to adapt ''on the fly'' the mesh from one time-step to another, and local time-stepping to retain the efficiency of explicit time integration even in the presence of local mesh refinement. 

\subsection{A posteriori error estimates}
\label{the:full-error analysis}
In \cite{GroteLakkisSantos:24:PosterioriErrorEstimatesWave} we derived the following a posteriori error estimates for the discretized time-dependent wave equation \eqref{eq:leapfrog} with mesh change:
\begin{equation}
\label{eq:errorboundenrg}
\max_{0 \leq n \leq N} \| U^n - u(.,t_n) \|_{\mathcal{A}} \leq \| \boldsymbol{e}(0) \|_{\text{erg}, \mathcal{A}}  +  C \left\{ 2 \sum_{m=1}^{2N} \zeta^m + \max_{1 \leq n \leq N} \varepsilon_0^n \right\},
\end{equation}
and
\begin{multline}
\label{eq:errorboundL2}
\max_{1 \leq n \leq N} \| V^{n-\frac{1}{2}} - v(.,t_{n-\frac{1}{2}})\|_{ L^2(\Omega)} 
\\
\leq
\| \boldsymbol{e}(0) \|_{\text{erg}, \mathcal{A}}  
+
C\left\{2 \sum_{m=1}^{2N} \zeta^m + \max_{1 \leq n \leq N} \varepsilon_1^n \right\}, 
\end{multline}
where the initial error, $\boldsymbol{e}(0)$, is defined as
$
 \boldsymbol{e}(0) := [u(0) - U^0, 
 v(t_{\frac{1}{2}}) - V^{\frac{1}{2}}]^{\top}
 $
 , and is a computable quantity.
Here the \textit{potential energy norm} is given by
\begin{equation}
    \Vert \phi \Vert_{\mathcal{A}} := \langle \mathcal{A} \phi \mid \phi \rangle^{\frac{1}{2}}
\end{equation}
and the \textit{wave-energy norm} by
\begin{equation}
\| \phi \|_{\mathrm{erg}, \mathcal{A}} := \langle \phi, \phi \rangle_{\mathrm{erg}, \mathcal{A}}^{1/2}.
\end{equation}
for the corresponding \textit{wave energy scalar product}
\begin{equation}
\langle \phi, \chi \rangle_{\mathrm{erg}, \mathcal{A}} := \langle \mathcal{A} \phi_0 \mid \chi_0 \rangle + \langle \phi_1, \chi_1 \rangle\text{ for }\phi = 
\begin{bmatrix}
\phi_0 \\ \phi_1
\end{bmatrix}, \, 
\chi = 
\begin{bmatrix}
\chi_0 \\ \chi_1
\end{bmatrix} \in H^1_{\Gamma_0} \times  L^2(\Omega).
\end{equation}

The upper bounds in \eqref{eq:errorboundenrg} and \eqref{eq:errorboundL2} involve standard a posteriori error estimators for elliptic problems with respect to the energy or the $L^2$ norm \cite{Verfurth:13:book:A-posteriori, AinsworthOden:00:book:A-posteriori} together with additional terms due to time integration or mesh change:

\vspace{1em}
\textbf{Elliptic error indicators:} (standard residual based error indicators)
\begin{align}\label{eq:elliptic_error_indicator}
\varepsilon_0^n &:= \mathscr{E}_{ \mathcal{A}} [U^n, \mathbbmss{V}_n], \\
\varepsilon_1^n &:= \mathscr{E}_{L^2} [V^{n-\frac{1}{2}}, \mathbbmss{V}_n] 
\end{align}
 with respect to the energy or the $L^2$-norm:
 \begin{align*}
\mathscr{E}_{\mathcal{A}}[w, \mathbbmss{V}]^2&:=
\sum_{K\in\mathcal{M}_{ \mathbbmss{W}}} \left\{ h_{\widehat K}^2\left\|A_{\mathbbmss{V}} w -\nabla \cdot(c(x) \nabla w)\big|_{K}\right\|^2_{\mathrm{L}_{2}(K)} \right. \\
&\mathrel{\phantom{=}} \left. + \frac{1}{2}h_{\widehat K}\left\|\llbracket c(x) \nabla w \big|_{K}\rrbracket\right\|^2_{\mathrm{L}_{2}\left(\partial K\right)} \right\} \\
\mathscr{E}_{L^2}[w, \mathbbmss{V}]^2&:=
\sum_{K\in\mathcal{M}_{ \mathbbmss{W}}} \left\{ h_{\widehat K}^4\left\|A_{\mathbbmss{V}} w -\nabla \cdot(c(x) \nabla w)\big|_{K}\right\|^2_{\mathrm{L}_{2}(K)}\right. \\
&\mathrel{\phantom{=}} \left.+\frac{1}{2}h_{\widehat K}^3\left\|\llbracket c(x) \nabla w \big|_{K}\rrbracket\right\|^2_{\mathrm{L}_{2}\left(\partial K\right)} \right\} 
\end{align*}
where $\widehat K\in \mathcal{M}_{ \mathbbmss{V}}$ denotes the smallest element 
which contains $K$.
Here $w$ typically denotes the (discrete) approximate solution 
to an elliptic problem while $\mathbbmss{V}$ denotes a finite-dimensional function space related to the problem. While in standard literature those
error estimates are stated on fixed meshes where the trial and test spaces
coincide, this is no longer the case in the presence of mesh change -- see Section 3.2 for further details.

In addition, the upper bounds in \eqref{eq:errorboundenrg} and \eqref{eq:errorboundL2} involve the following new error terms related to mesh change and the time discretization. 

\vspace{1em}
\textbf{Time accumulation indicators:}
\begin{equation}
\zeta^m := \int_{t_{\frac{m-1}{2}}}^{t_{\frac{m}{2}}} \left( (\mu_0^n + \vartheta_0^n(t))^2 + (\alpha^n + \mu_1^n + \delta^n(t) + \vartheta_1^n(t))^2 \right)^{\frac{1}{2}} \, \operatorname{d}\!t,
\end{equation}
for $n=\lfloor{\frac{m+1}{2}\rfloor}$ and $m = 1, \dots, 2N$, which result from 
the time-discretization, mesh change, local time-stepping and the numerical approximation of the right-hand side (forcing): 
\vspace{1em}
\textbf{Mesh-change indicators:} (nonzero only when the mesh changes)
\label{eq:def:mesh-change-indicator}
\begin{align*}
\mu_0^n &:= \left( \left\| \left[\Pi_{n} - \operatorname{Id}\right] U^{n-1} \right\|_{\mathcal{A}} + \mathscr{E}_{\mathcal{A}} \left[ [\Pi_{n} - \operatorname{Id}] U^{n-1}, \mathbbmss{V}_{n} \cap \mathbbmss{V}_{n+1} \right] \right) \Delta t^{-1}, \\
\mu_1^n &:= \left( \left\| [\Pi_{n+1} - \operatorname{Id}] V^{n-\frac{1}{2}} \right\|_{\operatorname{L}_2(\Omega)} + \mathscr{E}_{\operatorname{L}_2} \left[ [\Pi_{n+1} - \operatorname{Id}] V^{n-\frac{1}{2}}, \mathbbmss{V}_{n} \cap \mathbbmss{V}_{n+1}\right] \right) \Delta t^{-1} , \\
\mu_2^n &:= \left\| \left[\operatorname{Id} - \Pi_{n+1} \right] \tilde{A}_n U^n \right\|_{\operatorname{L}_2(\Omega)} + \mathscr{E}_{\operatorname{L}_2} \left[ \left[\operatorname{Id} - \Pi_{n+1} \right] \tilde{A}_n U^n, \mathbbmss{V}_{n+1}] \right];
\end{align*}
Note that the mesh-change indicators are typically small because only those elements
which change from one step to the next will contribute to them. In addition, we
only allow elements to change from one mesh to the next where those contributions are predictably small. 

\vspace{1em}
\textbf{LTS error indicators:} (due to using $\tilde{A}_n$ in scheme instead of $A_n$)
\begin{align}\label{eq:lts_error_indicator}
\alpha_0^n &:= \left\| \left[A_n - \tilde{A}_n \right] U^n \right\|_{\operatorname{L}^2(\Omega)},
\\
\alpha_1^n &:= \mathscr{E}_{L^2(\Omega)} \left[\tilde{A}_n U^n, \mathbbmss{V}_{n+1} \right],
\\
\alpha^n &:= \alpha_0^n + \alpha_1^n + \mu_2^n;
\end{align} 

\vspace{1em}
\textbf{Time-error indicators:} (mainly due to time discretization)
\begin{align}\label{eq:time_error_indicator_I}
\vartheta_0^n(t) &:= \Delta t^2 
\begin{cases}
     \left\Vert  \partial^2 V^{n-\frac{1}{2}} \frac{\ell_n(t)-1}{2} - \partial [A_{n-1} U^{n-1}]  q_{n-1}(t) \right\Vert_{\mathcal{A}} \\ \quad + \mathscr{E}_{\mathcal{A}} \left[ \partial^2 V^{n-\frac{1}{2}}\frac{\ell_n(t)-1}{2}, \mathbbmss{V}_{n-1} \cap\mathbbmss{V}_{n} \cap \mathbbmss{V}_{n+1} \right],\quad t \in I^{\prime}_{n-\frac{1}{2}}, \\[1em]
     \left\Vert \partial^2 V^{n-\frac{1}{2}} \frac{\ell_n(t)-1}{2} - \partial [A_{n} U^{n}]  q_{n}(t) \right\Vert_{\mathcal{A}} \\ \quad + \mathscr{E}_{\mathcal{A}} \left[ \partial^2 V^{n-\frac{1}{2}}\frac{\ell_n(t)-1}{2}, \mathbbmss{V}_{n-1} \cap\mathbbmss{V}_{n} \cap \mathbbmss{V}_{n+1} \right],\quad t \in I^{\prime}_{n},
\end{cases}
\end{align}
\begin{align}\label{eq:time_error_indicator_II}
\vartheta_1^n(t) &:= \Delta t^2 
\begin{cases}
     \left\Vert \frac{1}{2} \partial^2 U^n\ell_n(t)- \partial V^{n-\frac{1}{2}} q_{n-\frac{1}{2}}(t) \right\Vert_{ L^2(\Omega)},\qquad t \in I^{\prime}_{n}, \\[1em]
     \left\Vert \frac{1}{2} \partial^2 U^n\ell_n(t)- \partial V^{n-\frac{1}{2}} q_{n+\frac{1}{2}}(t) \right\Vert_{ L^2(\Omega)},\qquad t \in I^{\prime}_{n+\frac{1}{2}},
\end{cases}
\end{align}
where  $\ell_\nu(t)$ is the piecewise linear (in fact, affine) function in  $t$ satisfying 
\begin{equation}
\ell_\nu(t_\nu) = 1 \text{ and } \ell_\nu(t_\nu + k \Delta t) = 0, \qquad k \neq 0
\end{equation}
and the \textit{quadratic bubble} $q_{\nu}(t)$ is defined as the positive part
of the quadratic polynomial which vanishes at $t_{\nu \pm 1 / 2}$ and takes maximum $\frac{1}{8}$ at
$t_{\nu}$ :\\
\begin{equation*}
q_{\nu}(t):=\frac{\left(t-t_{\nu-1 / 2}\right)\left(t_{\nu+1 / 2}-t\right)}{2(\Delta t)^{2}} \mathbbm{1}_{\left[
\left|t-t_{\nu}\right|>\Delta t/2\right]} \quad \quad \nu=0,1 / 2,1, \dots, N-1 / 2, N. 
\end{equation*}
In addition, we denote
the centered difference in time at $t_\nu$ by
\begin{equation}
  \partial\phi^\nu
  :=
  \frac{\phi^{\nu+1} - \phi^{\nu-1}}{2 \Delta t}
\end{equation}
and
the centered second difference in time at $t_\nu$ by
\begin{equation}
  \partial^2\phi^\nu
  :=
  \frac{\phi^{\nu+1} - 2 \phi^\nu + \phi^{\nu-1}}{\Delta t^2}
\end{equation}

\vspace{1em}
\textbf{Data approximation indicator:} (due to a possibly nonzero source)
\begin{align}\label{eq:data}
\delta^n(t) := \| F^n - f(t) \|_{ L^2(\Omega)}.
\end{align}

\subsection{Compatible mesh change}\label{sec:compatible_mesh_change}
Mesh change during any time-dependent simulation cannot be arbitrary
without risking a significant loss in accuracy \cite{Dupont:82:article:Mesh}.
Indeed both mesh coarsening or refinement will generally increase the numerical error
when transferring the FE solution via interpolation or projection from one mesh to another. To ensure that local mesh refinement never leads to a loss in accuracy, 
we restrict mesh change to compatible meshes only. Moreover, we shall 
allow coarsening only where appropriate to minimize the resulting inherent information loss.

Hence we assume that the domain $\Omega$ is a polytope and that it can be partitioned into simplices exactly with the coarsest mesh, $\mathcal{M}_0$, called the \emph{macro triangulation} where every element of $\mathcal{M}_0$ is "ready" to be bisected (following the newest vertex bisection algorithm in 2-D and the Kossaczký algorithm in 3-D). We call two meshes $\mathcal{M}_1$  and $\mathcal{M}_2$ compatible, if every element $K\in \mathcal{M}_1$ is either an element of $\mathcal{M}_2$ or a union of elements thereof, and vice-versa. 
Thus any two meshes $\mathcal{M}_{\mathbbmss{V}}$ and $\mathcal{M}_{\mathbbmss{W}}$
with their corresponding FE spaces ${\mathbbmss{V}}$ and ${\mathbbmss{W}}$ 
that were obtained via refinement by bisection of $\mathcal{M}_0$ are compatible. 
Moreover, the two FE subspaces ${\mathbbmss{V}} + {\mathbbmss{W}}$ and 
${\mathbbmss{V}} \cap {\mathbbmss{W}}$ correspond to the FE spaces
associated with the coarsest common refinement and finest common coarsening, respectively, of $\mathcal{M}_{\mathbbmss{V}}$ and $\mathcal{M}_{\mathbbmss{W}}$.

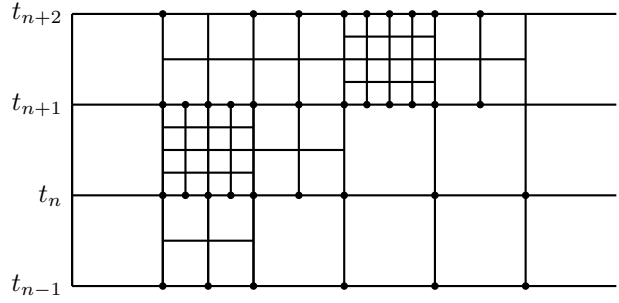
\begin{figure}[t]
\begin{center}
\begin{tikzpicture}[scale=1.2]

\foreach \y/\label in {0/{$t_{n-1}$},1/{$t_n$},2/{$t_{n+1}$},3/{$t_{n+2}$}} {
  \draw[thick] (0,\y) -- (6,\y);
  \node[left] at (0,\y) {\label};
}

\foreach \x in {0,1,2,3,4,5} {
  \draw[thick] (\x,0) -- (\x,1);
}

\foreach \x in {0,1,1.5,2} {
  \draw[thick] (\x,0) -- (\x,2);
}

\draw[thick] (1.5,0) -- (1.5,1);
\draw[thick] (1,0.5) -- (2,0.5);

\foreach \x in {0,1,1.25,1.5,1.75,2,2.5,3,4,5} {
  \draw[thick] (\x,1) -- (\x,2);
}

\draw[thick] (1,1.5) -- (3,1.5);

\draw[thick] (1,1.25) -- (2,1.25);

\draw[thick] (1,1.75) -- (2,1.75);

\foreach \x in {0,1,1.5,2,2.5,3,3.25,3.5,3.75,4,4.5,5} {
  \draw[thick] (\x,2) -- (\x,3);
}

\draw[thick] (1,2.5) -- (5,2.5);

\draw[thick] (3,2.25) -- (4,2.25);

\draw[thick] (3,2.75) -- (4,2.75);

\foreach \x/\y in {
  1/0,1.5/0,2/0,3/0,4/0,5/0,
  1/1,1.25/1,1.5/1,1.75/1,2/1,2.5/1,3/1,4/1,5/1,
  1/2,1.25/2,1.5/2,1.75/2,2/2,2.5/2,3/2,3.25/2,3.5/2,3.75/2,4/2,4.5/2,
  1/3,2/3,2.5/3,3/3,3.25/3,3.5/3,3.75/3,4/3,4.5/3
} {
  \filldraw[black] (\x,\y) circle (1pt);
}

\end{tikzpicture}
\end{center}
 \caption{Space-time refined compatible meshes.}
    \label{fig:refined_grid}
\end{figure}

At time-step $t_n \mapsto t_{n+1}$, any element $K$ in the underlying FE mesh 
$\mathcal{M}_n$ will either remain as is, or change through local coarsening or refinement, as shown in Fig. \ref{fig:refined_grid}. If $\widehat K\in\mathcal{M}_n$ is refined by bisection, 
$\widehat K =  K_{-} \cup K_{+}$, we simply interpolate the FE approximation on the two new elements $K_{-}, K_{+} \in \mathcal{M}_{n+1}$. Hence no additional error results from compatible local mesh refinement. On the other hand, if two neighboring elements $K_{-}, K_{+} \in \mathcal{M}_n$, both children of the same coarser parent element $\widehat K$ in the refinement tree, are marked for coarsening, they will be replaced by $\widehat K \in\mathcal{M}_{n+1}$. Here to minimize the inherent information loss, we estimate in advance the potential loss in accuracy due to coarsening by computing the following {\em coarsening pre-indicators} $\beta_K$. 

Let $y(x) \in {\mathbbmss{V}}_n$ be a current FE approximation associated with mesh $\mathcal{M}_n$
and assume that $K_{-}, K_{+} \in \mathcal{M}_n$ are both children of the same coarser parent element $\widehat K \in\mathcal{M}_{n+1} $, with $\widehat K = K_- \cup K_+$. Next, let $\hat y(x) = \Pi_{n+1} y(x)$ denote its FE approximation in  ${\mathbbmss{V}}_{n+1}$  obtained via interpolation (or projection). Then we define for each element
$K \subset \widehat K$ marked for possible coarsening 
the {\em coarsening pre-indicators}
\begin{eqnarray}\label{eq:beta}
    \beta_K^0 &:=& \left\| \left[\Pi_{n+1} - \operatorname{Id}\right] y \right\|_{\mathcal{A}, K} =  \left\| (y -\hat{y})|_K\right\|_{\mathcal{A}} \\
      \beta_K^1 &:=& \left\| \left[\Pi_{n+1} - \operatorname{Id}\right] y \right\|_{L^2(K)} =  \left\| (y -\hat{y})|_K\right\|_{L^2}.
\end{eqnarray}
Owing to the compatibility of the two meshes $\mathcal{M}_n$ and $\mathcal{M}_{n+1}$,
both pre-indicators are easily computed for each $K_-, K_+ \subset \widehat K$ as
each FE basis function $\Phi_m^{n+1} \in {\mathbbmss{V}}_{n+1}$ restricted to $\widehat K$ is a linear combination of FE basis functions $\Phi_m^n \in {\mathbbmss{V}}_{n}$. Moreover, for a Lagrangian FE basis, $y(x) - \hat y (x)$ vanishes at all common nodes $x_j$ and thus reduces to a linear combination of FE basis functions over the remaining nodes. 

Following \cite{LakkisPryer:12:article:Gradient} we now
exemplify this point by considering the simplest situation of a one-dimensional
FE approximation with continuous piecewise linear polynomials. Again let
$K_- = (x_{m-1}, x_m)$, $K_+ = (x_{m}, x_{m+1})$ both marked for coarsening 
and thus possibly replaced by the single element $\widehat K = (x_{m-1}, x_{m+1})$.
For $x\in\widehat K$ we have
\begin{eqnarray*}
    y(x) &=& y_{m-1}\Phi^n_{m-1}(x) + y_{m}\Phi^n_{m}(x) + y_{m+1}\Phi^n_{m+1}(x) \\
    \hat y(x) &=& y_{m-1}\Phi^{n+1}_{m-1}(x) + y_{m+1}\Phi^{n+1}_{m+1}(x), 
\end{eqnarray*}
where $\Phi_m^n(x) \in {\mathbbmss{V}}_{n}$ denotes a standard ''hat-function''
with $\Phi_m^n(x_j) = \delta_{m,j}$. Next, we note that the ''coarser'' basis
functions are linear combinations of the ''finer'' ones
for $x\in\widehat K$,
\begin{eqnarray*}
    \Phi^{n+1}_{m-1}(x) &=& \Phi^n_{m-1}(x) + \frac 12\Phi^n_{m}(x) \\
      \Phi^{n+1}_{m+1}(x) &=& \Phi^n_{m+1}(x) + \frac 12\Phi^n_{m}(x),
\end{eqnarray*}
where we have used that the two children elements $K_-, K_+$, previously obtained via bisection of $\widehat K$, are equally sized.
We thus easily calculate
\begin{equation}
    y(x) - \hat y(x) = \left[y_m - \frac 12 (y_{m-1}+y_{m+1}) \right] \Phi_m^n(x), \qquad x\in \widehat K,
\end{equation}
which immediately yields the pre-indicators
\begin{eqnarray*}
    \beta_{K_-}^1 &=& \left| y_m - \frac 12 (y_{m-1}+y_{m+1})\right| \left\| \Phi^n_m\right\|_{L^2{(K_-)}} \\
    \beta_{K_+}^1 &=& \left| y_m - \frac 12 (y_{m-1}+y_{m+1})\right| \left\| \Phi^n_m\right\|_{L^2{(K_+)}}, 
\end{eqnarray*}
and similarly for $\beta_{K_-}^0$ and $\beta_{K_+}^0$ with the $L^2$-norm replaced by the energy-norm. 

When the pre-indicators of two neighboring elements $K_{-}$ and $K_+$ are sufficiently small, those elements are subsequently replaced by a single coarser element $\widehat K$. If standard polynomial interpolation
is used for $\Pi_{n+1}$ in the newly coarsened FE space ${\mathbbmss{V}}_{n+1}$, its remaining nodal values will not change.
 By allowing for coarsening only in those elements where the pre-indicators are
small, we thus also avoid any detrimental drastic increase in the mesh-change
indicators $\mu_0^n, \mu_1^n$ and $\mu_2^n$ from Section 3.1 .

\subsection{Local time-stepping}\label{sec:lts}
Here we briefly recall the leapfrog (LF) based local time-stepping (LF-LTS) method used for time integration during the adaptive algorithm to circumvent the bottleneck otherwise caused by local refinement on standard explicit 
time-stepping methods. 
The original LF-LTS method for the
numerical solution of the second-order wave equations was
proposed for homogeneous right-hand sides in \cite{DiazGrote:09:article:Energy-Conserving} and for inhomogeneous
right-hand sides in \cite[Sect. 4.1]{GroteMitkova:10:article:Explicit-Local}.
Inside the locally refined region, the latter applies standard
LF time-marching with a smaller time-step, which also implies evaluating the inhomogeneous right-hand side at all intermediate times.
Although we focus on second-order LTS-LF here, we note that 
higher-order versions are available \cite{DiazGrote:09:article:Energy-Conserving}.
For $f=0$ they also conserve (a discrete version of) the energy. 
Optimal convergence rates for the
LF-LTS method from \cite{DiazGrote:09:article:Energy-Conserving} with $p$ local time steps were derived for a
conforming FEM discretization, albeit under a sub-optimal CFL 
condition where $\Delta t$
in fact depends on the smallest elements in the mesh \cite{GroteMehlinSauter:18:article:Convergence}. 

To prove optimal $L^{2}$ convergence rates under a CFL condition independent of
$p$, a stabilized algorithm LF-LTS$(\nu)$ was recently introduced in
\cite{GroteMichelSauter:21:article:Stabilized}, and also independently in \cite{CarleHochbruck:22:article:Error-Analysis}. 
Here, $\nu \geq 0$ denotes a small stabilization parameter; typically, we set $\nu=0.01$. 
Stability and convergence for the stabilized version LF-LTS$(\nu)$
were proved for homogeneous right-hand sides in
\cite{GroteMichelSauter:21:article:Stabilized} under a CFL condition independent of the coarse/fine mesh ratio. 

Stabilization was also recently
introduced into the original LTS-LF$(0)$ algorithm from  \cite{GroteMitkova:10:article:Explicit-Local} for nonzero source terms, which again led to 
optimal $L^2$-convergence rates under a CFL condition independent of
the coarse-to-fine mesh ratio \cite{GroteMichelSauter:24:ExplicitLocalTimeSteppingInhomogeneous}.  
Similar optimal $L^2$-convergence rates for a somewhat simpler "split-LFC" variant, which omits intermediate source evaluations inside the refined region, were
recently proved in \cite{CarleHochbruck:24:article:Error}.

In \eqref{eq:leapfrog}, $\widetilde{A_{n}}$ denotes the perturbed bilinear form introduced by local time-stepping, which computes the solution 
at $t_{n+1}$ by using a
smaller time-step $\Delta\tau=\Delta t/p$ inside the regions of local
refinement; here, $p\geq2$ denotes the \textquotedblleft
coarse\textquotedblright\ to \textquotedblleft fine\textquotedblright%
\ time-step ratio. It is given by
\begin{equation}
\label{eq:LTS_A_nostab}
\widetilde{A_{n}}:=A_{n} P_p(\Delta t^2\Pi_{n}^{\mathrm{f}} A_{n}),
\end{equation}
where the polynomial $P_{p}(x)$ is given by
\[
P_{p}(x) = \frac{2}{x} \left( 1 - T_p \left( 1 - \frac{x}{2 \,p^2} \right) \right) ,
\]
with $T_p(x)$ the Chebyshev polynomial of order $p$, recursively defined as
\[
T_p(x) = 2 \, x \,T_{p-1}(x) - T_{p-2}(x), \quad p\geq 2, \qquad  T_0(x) = 1, \quad T_1(x) = x.
\]
For $p=2$, for instance, $T_2(x) = 2x^2-1$ and thus $P_2(x) = 1 - x/16$,
so that \eqref{eq:LTS_A_nostab} indeed reduces to \eqref{eq:LTS_A}.

In practice we never explicitly compute the polynomial $P_p(x)$ but
instead use the "leapfrog-like" three-term recursion satisfied by
the Chebyshev polynomials. The $p$ sub-iterations then correspond to 
$p$ local time-steps with $\Delta \tau := \Delta t/p \sim h_{\min}$ 
that affect only the unknowns in refined region, as those additional $p$ multiplications with $\Pi_{n}^{\mathrm{f}} A_{n}$ only modify unknowns 
inside the refined part. 
For $p = 1$, $P_1(x) = 1$ and the method reduces to the standard leapfrog scheme.

Since instabilities are possible for certain values of  $\Delta t$, 
we generally replace $P_p(x)$ in \eqref{eq:LTS_A_nostab} by
so-called "damped" Chebyshev polynomials
\[
P_{p,\nu}(x) = \frac{2}{x} \left( 1 - \frac{T_p \left( \delta_{p,\nu} - \frac{x}{\omega_{p,\nu}} \right)}{T_p(\delta_{p,\nu})} \right),
\]
where
\[\nu >0,\quad \delta_{p,\nu} = 1 + \frac{\nu}{p^2}, \quad
\omega_{p,\nu} = 2 \, \frac{T_p^{\prime}(\delta_{p,\nu})}
{T_p(\delta_{p,\nu})} .\]
The stabilized version is denoted by LF-LTS$(\nu)$;
for $\nu = 0$, the original LF-LTS$(0)$ method is recovered.

 \subsection{Space-time adaptive algorithm}
\label{sect:algo}
Given a fixed ''coarse'' FE mesh, $\mathcal{M}_0$, which discretizes the computational domain $\Omega$, we first choose
a constant global time-step $\Delta t$ which satisfies the CFL stability condition imposed by $\mathcal{M}_0$ on the (standard) leapfrog method.
Next, the mesh is adapted to the initial conditions $u_0, v_0$ - see \cref{alg:initial} -- which yields the initial mesh $\mathcal{M}_{1}$ used for time integration. Now, during each time step, we proceed as follows. First, we tentatively set $\mathcal{M}_{n+1}$ to the current mesh $\mathcal{M}_n$ and compute $U^{n+1}$ using the LTS-LF method inside locally refined regions and the standard LF method elsewhere. Next, 
we evaluate the elliptic error indicators $\varepsilon_0$ and $\varepsilon_1$ in \eqref{eq:elliptic_error_indicator} from $U^{n+1}$ and refine $\mathcal{M}_{n+1}$ accordingly, if needed. We then transfer $U^{n-1}$ and $U^n$ to the new finite element space $\mathbbmss{V}_{n+1}$, associated with $\mathcal{M}_{n+1}$ via interpolation or $\operatorname{L}^2$-projection, and repeat the computation of $U^{n+1}$. This process -- refining  $\mathcal{M}_{n+1}$ and recomputing $U^{n+1}$ -- is iterated until $\varepsilon_0$ or $\varepsilon_1$ falls below the threshold 
$\tol{H}/{N}$. At that point, $U^{n+1}$ is accepted on the current mesh $\mathcal{M}_{n+1}$. Now, the algorithm enters a coarsening step: for each $K\in\mathcal{M}_{n+1}$, we evaluate the pre-indicators $\beta_K$ from \eqref{eq:beta} and mark those elements with  $\beta_K\leq\tol{C}$. After coarsening, both $U^n$ and $U^{n+1}$ are transfered to the new, coarsened space $\mathbbmss{V}_{n+1}$. The entire procedure is repeated until the final time $T$ is reached.

\begin{algorithm}[H]
\caption{Space-time Adaptive LF-FEM}\label{alg:adaptive}
\begin{algorithmic}[1]
\Require initial conditions $u_0, v_0$, initial mesh $\mathcal{M}_{0}$, $\tol{H}$, $\tol{C}$, final time $T$, number of time-steps $N$, D\"orfler marking threshold $\theta > 0$
\Ensure adapted meshes $(\mathcal{M}_n)_{n=0, \ldots, N}$ and discrete solution $(U^n)_{n=0, \ldots, N}$ such that $\Vert U - u \Vert < \tol{}$
\Procedure{adapt}{$T, \tol{H}, \tol{C}, \Delta t, \theta, u^0, v^0$} 
\State $(U^0, U^1, \mathcal{M}_1) \gets \Call{initialize}{u^0,v^0, \mathcal{M}_0, \theta}$
\For{$n = 1:N-1$} 
 \State $\mathcal{M}_{n+1} \gets \mathcal{M}_{n} $
 \State $\mathcal{R} \gets \emptyset$
 \State set $ 
 \varepsilon_i^{n+1}, \mu_j^{n+1},\alpha_i^{n+1}, \vartheta_i^{n+1}, \delta^{n+1}
    \gets \nicefrac{2\tol{H}}{N}$ \Comment{$
    i=0,1, j=0,1,2
    $}
  \While{$
\begin{aligned}
\min (
\varepsilon_{i}^{n+1},\mu_j^{n+1}, \alpha_i^{n+1}, \vartheta_i^{n+1}, \delta^{n+1}
)
     > \nicefrac{\tol{H}}{N}
\end{aligned}
$}
        \State compute $U^{n+1} \gets 2U^n - U^{n-1} + ( F^n
- \tilde{A}_n U^n) \Delta t^2$
            \State compute $
            \varepsilon_i^{n+1}, \mu_j^{n+1}, \alpha_i^{n+1},\vartheta_i^{n+1}, \delta^{n+1}
$ \Comment{store $(\eta_i^K)_K, i=0,1$}

            \State $\mathcal{D} \gets \emptyset$
            \For{$i=0,1$} 
            \State $s \gets 0$
            \For{$K \in \mathcal{M}_{n+1}$} \\
            \Comment{where $\mathcal{M}_{n+1}$ is in decreasing order of $\eta_i^K$ }
         \If{$s \leq \theta \, ((
         \varepsilon_i^{n+1})^2 + (\mu_j^{n+1})^2 +  (\alpha_i^{n+1})^2 + (\vartheta_i^{n+1})^2 + (\delta^{n+1})^2
         )$}
         \State $s \gets s + ( \eta_i^K)^2$
         \State add $K$ to $\mathcal{D}$
         \EndIf
         \EndFor
         \EndFor
         \State $(\mathcal{M}_{n+1},\mathcal{D}) \gets \Call{refine}{\mathcal{M}_{n+1},\mathcal{D}} $
         \State add $\mathcal{D}$ to $\mathcal{R}$
    \EndWhile
    \State $\mathcal{C} \gets \emptyset$ 
    \For{$K \in \mathcal{M}_{n+1}\setminus \{ \mathcal{R} \}$}
    \If{$ \max{(\beta_K^0,\beta_K^1)} <\tol{c}$}
    \State add $K \; \text{to} \; \mathcal{C}$ 
    \EndIf
    \EndFor
     \State $(\mathcal{M}_{n+1},\mathcal{C}) \gets \Call{coarsen}{\mathcal{M}_{n+1},\mathcal{C}} $
\EndFor
\EndProcedure
\end{algorithmic}
\end{algorithm}
\begin{remarks}
    \begin{enumerate}
        \item In line 2 the initialization function is used to determine the initial mesh based on the initial conditions - see \cref{alg:initial}.
        \item The quantities $\varepsilon_0$ and $\varepsilon_1$ denote the elliptic error indicators in \eqref{eq:elliptic_error_indicator}.
        \item In line 15, we use a Dörfler marking strategy for refinement \cite{Verfurth:13:book:A-posteriori}.
        \item The set $\mathcal{D}$ contains the elements to be refined, whereas the set $\mathcal{C}$ contains those elements to be coarsened.
        \item We ensure within each time step (for-loop in line 3) that we never coarsen any element just recently refined.  
        \item After computing the new mesh $\mathcal{M}_{n+1}$, the previous solutions $U^n$ and $U^{n-1}$ are transfered to $\mathbbmss{V}_{n+1}$.
        \item 
        In line 8, the leapfrog based explicit local time-stepping method from Section 3.3 is used with a time-step ratio $p$ given by the local mesh size ratio. Thus for any refined region we use the value of $p$ determined by the ratio of the coarse to the smallest mesh size in that subregion; note that $p$
	may vary from one refined subregion to another.  In situations where
        the refined region itself contained yet another sub-region of much smaller mesh-size, a hierarchical multi-level approach could be used instead for even higher efficiency
        \cite{DiazGrote:15:article:Multi-Level-Explicit,Rietmann:17:NewmarkLocalTime}.
    \end{enumerate}
\end{remarks}
Before starting the actual time integration, the initial uniform mesh $\mathcal{M}_0$ is adaptively refined to better resolve the initial conditions \( u_0 \) and \( v_0 \), resulting in the mesh \( \mathcal{M}_1 \). First, we estimate the local error indicator \( \eta_0 \) for \( u_0 \) by computing the element-wise product of the local mesh size and the \( \operatorname{L}^2 \)-norm of its second spatial derivative.
Elements for which \( \eta_0(K) > \theta \max_K \{\eta_0(K)\} \) are marked for refinement and refined using a Dörfler marking strategy \cite{Verfurth:13:book:A-posteriori}. The same procedure is applied to \( v_0 \) using a second error indicator \( \eta_1 \), which eventually yields the initial mesh \( \mathcal{M}_1 \). The initial data \( u_0 \) and \( v_0 \) are then discretized on \( \mathcal{M}_1 \) via interpolation. Finally, the first time step \( U^1 \) is computed using a Taylor expansion.

\begin{algorithm}[H]
\caption{Initialize}\label{alg:initial}
\begin{algorithmic}[1]
\Require $u^0,v^0,\mathcal{M}_0$,$\theta$
\Ensure $U^0, U^1, \mathcal{M}_1$
\Procedure{initialize}{$u^0,v^0,\mathcal{M}_0$,$\theta$}
\State $\boldsymbol{\eta}^0 \gets \Call{estimate}{\mathcal{M}_0,h_K \Vert \partial_{xx} u^0 \Vert_{L_2(K)}}$
\State $\mathcal{D} \gets \Call{mark}{\mathcal{M}_0, \theta, \boldsymbol{\eta}^0}$
        \State $(\mathcal{M}_{0},\mathcal{D}) \gets \Call{refine}{\mathcal{M}_{0},\mathcal{D}} $
        \State $\boldsymbol{\eta}^1 \gets \Call{estimate}{\mathcal{M}_0,h_K \Vert \partial_{xx} v^0 \Vert_{L_2(K)}}$
\State $\mathcal{D} \gets \Call{mark}{\mathcal{M}_0, \theta, \boldsymbol{\eta}^1}$
        \State $(\mathcal{M}_{1},\mathcal{D}) \gets \Call{refine}{\mathcal{M}_{0},\mathcal{D}} $
        \State $U^0 \gets \Call{discretize}{u^0, \mathcal{M}_1}$
        \State $V^0 \gets \Call{discretize}{v^0, \mathcal{M}_1}$
        \State $U^1 \gets U^0 + \Delta t V^0 + \frac{\Delta t^2}{2}(F^0-A^0U^0)$
\EndProcedure
\end{algorithmic}
\end{algorithm}
\section{Numerical results}
Here we present a series of numerical experiments that confirm the optimal convergence rates of the LF--LTS--FEM method \eqref{eq:leapfrog} with mesh change and demonstrate the effectiveness of the space--time adaptive algorithm from Section~\ref{sect:algo}. 
First, we consider the solution of the inhomogeneous wave equation \eqref{eq:waveeq} on a predefined, time-evolving mesh under successive mesh refinement to verify convergence. 
Next, we apply the space--time adaptive LF--FEM algorithm to compute two one-dimensional test cases: a right-moving Gaussian pulse and a Gaussian pulse that splits into two symmetric waves traveling in opposite directions. 
Finally, we consider a two-dimensional example, where 
a Gaussian pulse initially centered at $(0.4,0.6)$ propagates across an L-shaped domain.

In all one-dimensional experiments, we solve the wave equation \eqref{eq:waveeq} in $\Omega = (-10,10)$ with homogeneous Dirichlet boundary conditions, i.e.\ $\Gamma = \Gamma_D$, and wave speed $c \equiv 1$. 
We use standard piecewise linear $H^1$-conforming finite elements on nonuniform meshes with mass-lumping in space and the leapfrog-based local time-stepping (LF-LTS) method with global time step $\Delta t$ (with added stabilization $\nu = 0.01$) from Section~\ref{sec:lts}. 

Since the entire time-stepping procedure is fully explicit---no linear systems are ever solved---the computational complexity scales linearly with the total number of degrees of freedom in both space and time.
For the one-dimensional adaptive LF--FEM experiments we choose the parameters $\tol H = 20$, $\tol C = 10^{-4}$, and $\theta = 0.8$
whereas in the two-dimensional experiment in Section~\ref{sec:pulse} we 
set $\tol C = 0.01$. 

\subsection{Forced wave}
To verify the convergence of the LF-LTS-FEM method \eqref{eq:leapfrog}, we first apply it to the inhomogeneous wave equation \eqref{eq:waveeq} on a sequence of predefined but time-varying meshes which follow a right-moving wave.
The nonzero source $f(x,t)$ and the initial conditions $u_0, v_0$ are set to match the linearly increasing right-moving Gaussian pulse, 
\begin{equation*}
  u(x,t) = t \, e^{-4(x-1-t)^2},
\end{equation*}
initially centered about $x = 1$.
  \par
\begin{figure}[t]
  \begin{subfigure}[t]{.5\textwidth}
    \centering
    \includegraphics[width=1\linewidth,height=4.3cm]{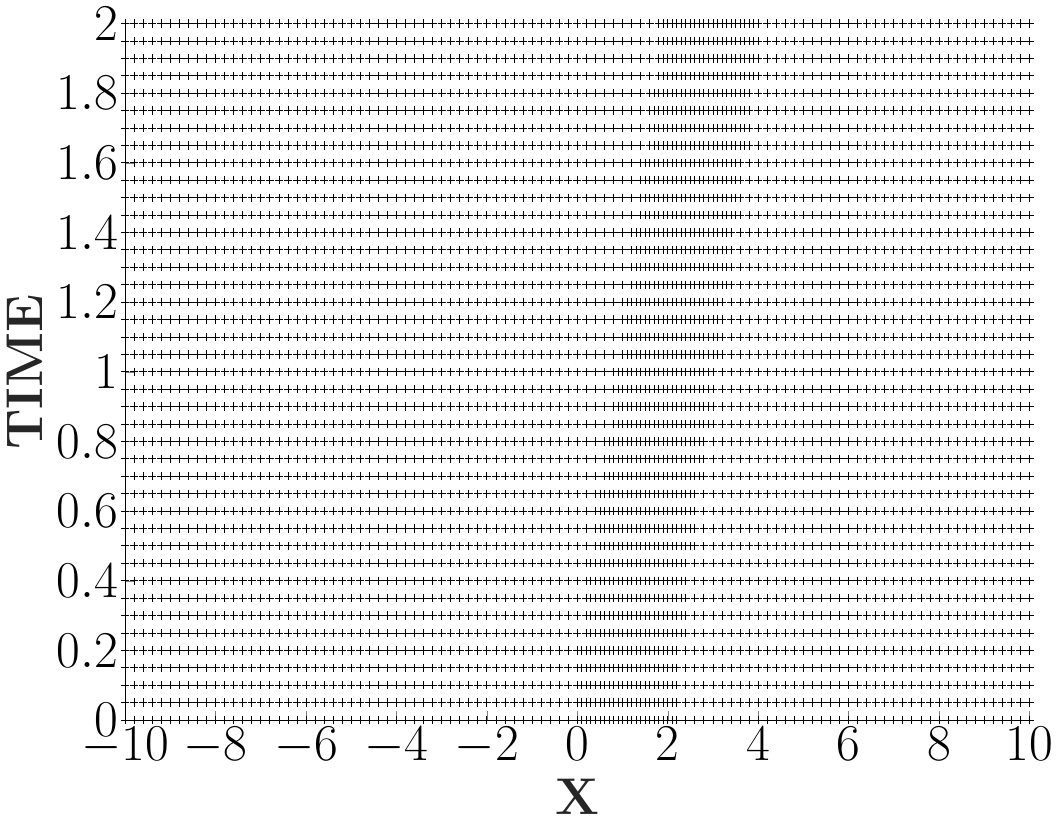}  
    \caption{Forced wave: Time-evolving mesh.}
    \label{fig:mesh_fixed}
  \end{subfigure}
  \hfill
  \begin{subfigure}[t]{.49\textwidth}
    \centering
    \includegraphics[width=1\linewidth]{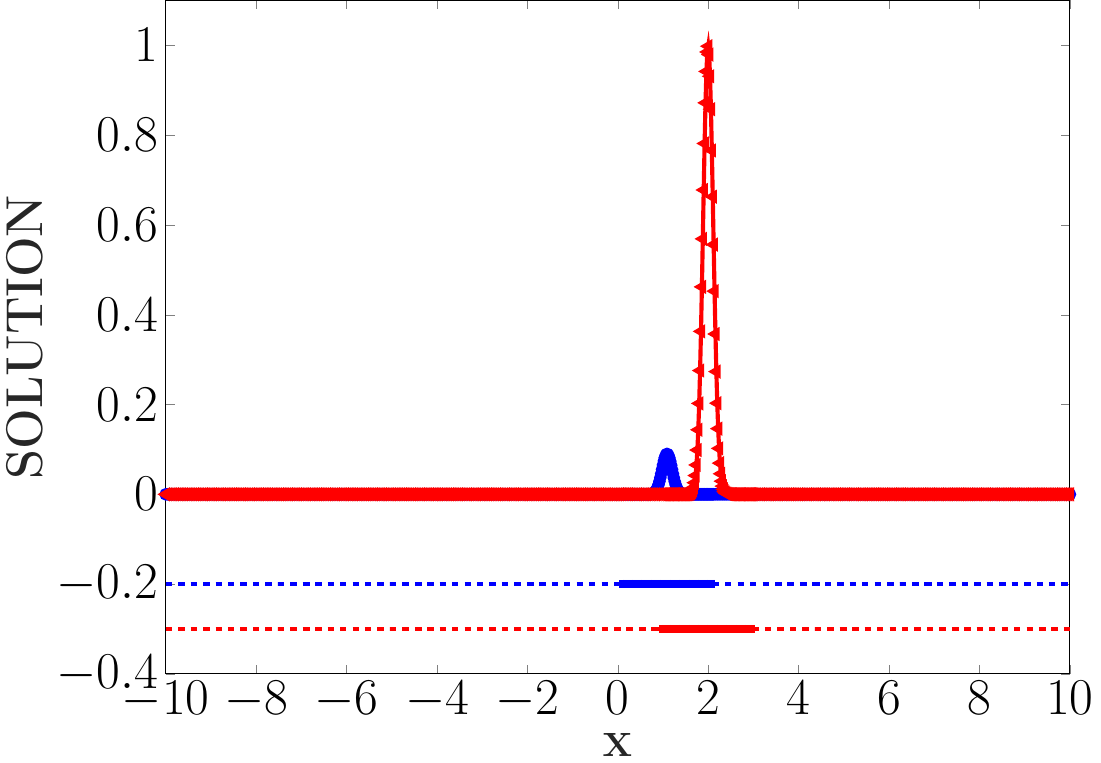}  
    \caption{Numerical solution and refined mesh at time $t=0.04$ (blue) and $t = 1$ (red).}
    \label{fig:solution_fixed}
  \end{subfigure}
\end{figure}

  At every discrete time $t_n$, the FE mesh $\mathcal{M}_n$ separates into a coarse part, $\mathcal{M}_n^c$, of
  constant mesh-size $h^c=h$, and a fine part, $\mathcal{M}_n^f$, of
  constant mesh-size $h^f=h^c / 2$; here, $h=h^c$ and $h^f$ themselves do not depend on time. The coarse and refined parts of the initial mesh, $\mathcal{M}_1$,
  correspond to 
  $\Omega^c_1 =[-10, 0] \cup [2,10]$ and $\Omega^f_1 = [0, 2]$, respectively. 
    The refined part of the mesh, $\mathcal{M}_n^f$, ``follows'' the Gaussian
  pulse propagating rightward with speed one across $\Omega$, as the mesh together with the associated FE space $\mathbbmss{V}_n$ change
whenever the
  elapsed time from the previous mesh change is greater than
  $h^c$. The resulting space-time mesh
  is shown in Fig. \ref{fig:mesh_fixed}.
  To ensure stability, we let the LF-LTS method take two local time-steps of
  size $\Delta t/2$ inside $\Omega^f_n$ during each global time-step of size $\Delta t$
  inside $\Omega^c_n$.

  During mesh change, two subsequent meshes $\mathcal{M}_n$ and
  $\mathcal{M}_{n+1}$ always remain compatible -- see Section \ref{sec:compatible_mesh_change}; hence,
  no additional discretization error occurs inside new elements from
  refinement. During coarsening, however, the removal of the common node at the
  interface between two fine elements introduces an additional discretization error, 
  which is kept small by allowing coarsening only where it is nearly zero. 
  The global time-step $\Delta$ corresponds to the CFL stability limit of
  a uniform mesh with identical mesh-size $h = h^c$.

\begin{figure}[t]
\begin{subfigure}[t]{.475\textwidth}
  \centering
  \includegraphics[width=1\linewidth]{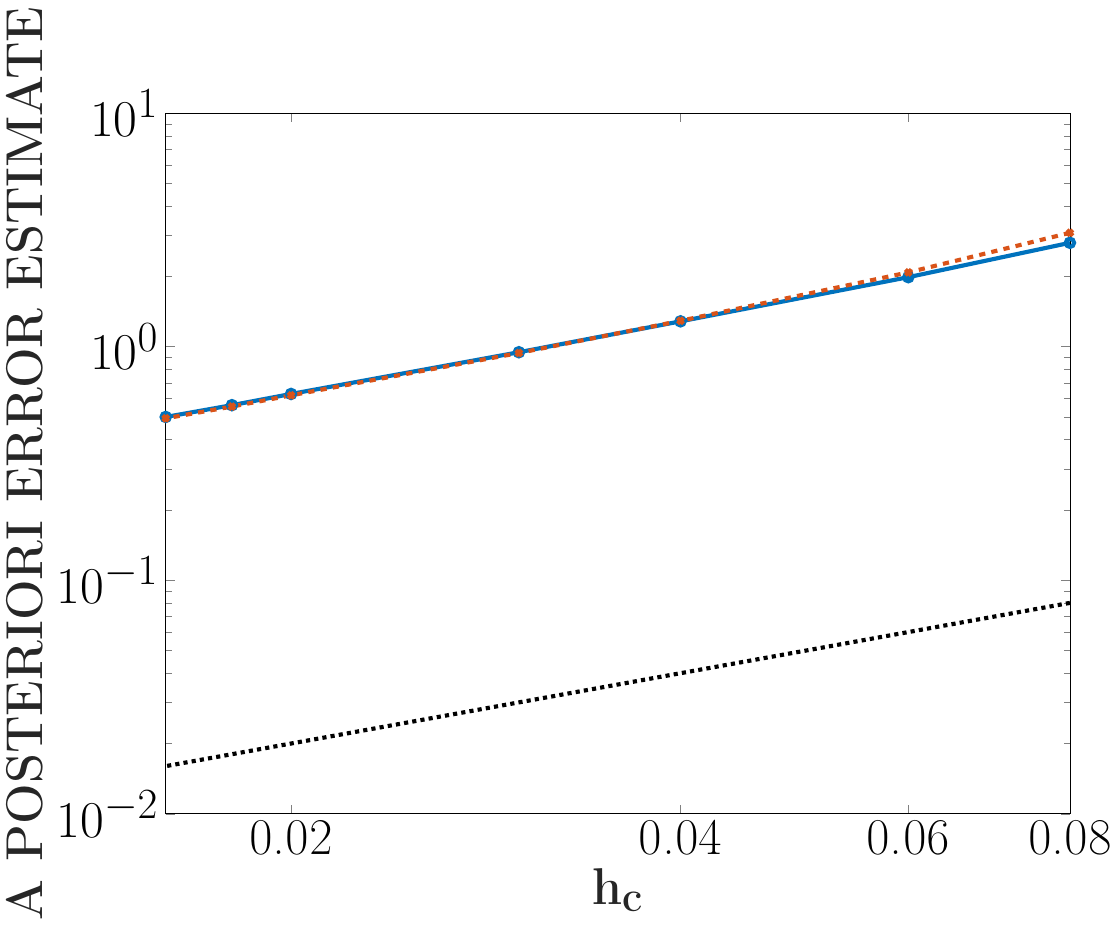}  
  \caption{Forced wave: A posteriori error estimates
\eqref{eq:errorboundenrg} (solid line with circles),
\eqref{eq:errorboundL2} (dashed line with squares), and
reference convergence rate $\mathcal{O}(h)$ (dotted line).}
  \label{fig:aposteriori_fixed}
\end{subfigure}
\hfill
\begin{subfigure}[t]{.475\textwidth}
\centering
  \includegraphics[width=1\linewidth]{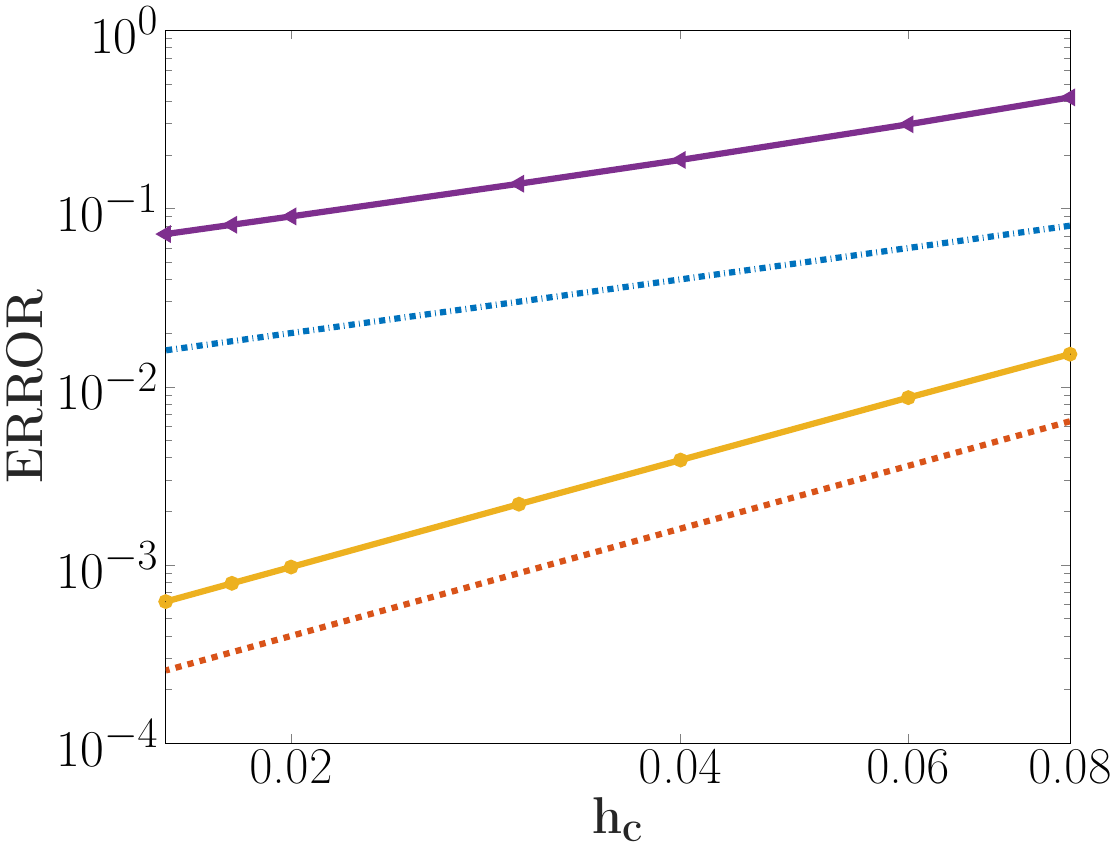}  
  \caption{Relative error in the
    energy-norm (solid line with diamonds) and the $\operatorname{L}_2(\Omega)$-norm (solid line with circles), 
    rates $\mathcal{O}(h)$ (blue dash-dot) and $\mathcal{O}(h^2)$ (red dash-dot).}
  \label{fig:convergence_fixed}
\end{subfigure}
\newline
\begin{subfigure}[t]{.475\textwidth}
  \centering
  \includegraphics[width=1\linewidth]{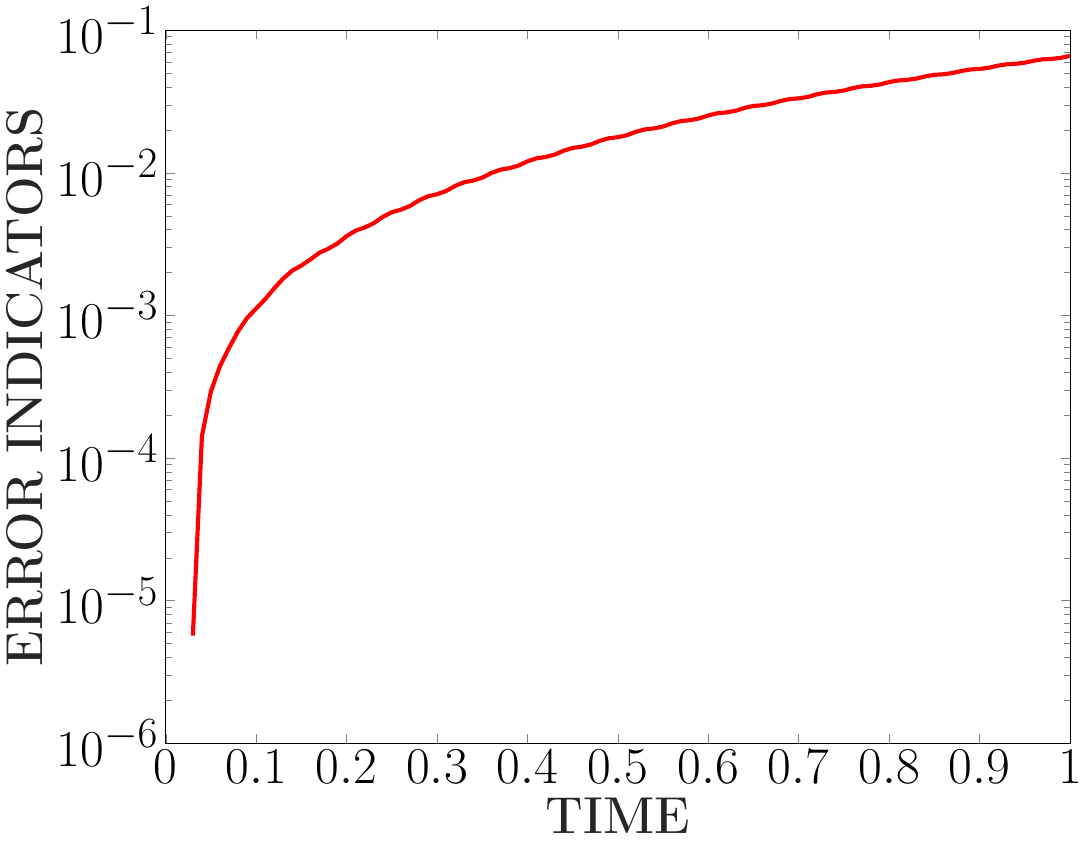}  
 \caption{Error indicator $\vartheta_0^n$ in \eqref{eq:time_error_indicator_I} vs. time.}
  \label{fig:r0-terms_accu_fixed}
\end{subfigure}
\hfill
\begin{subfigure}[t]{.475\textwidth}
  \centering
  \includegraphics[width=1\linewidth]{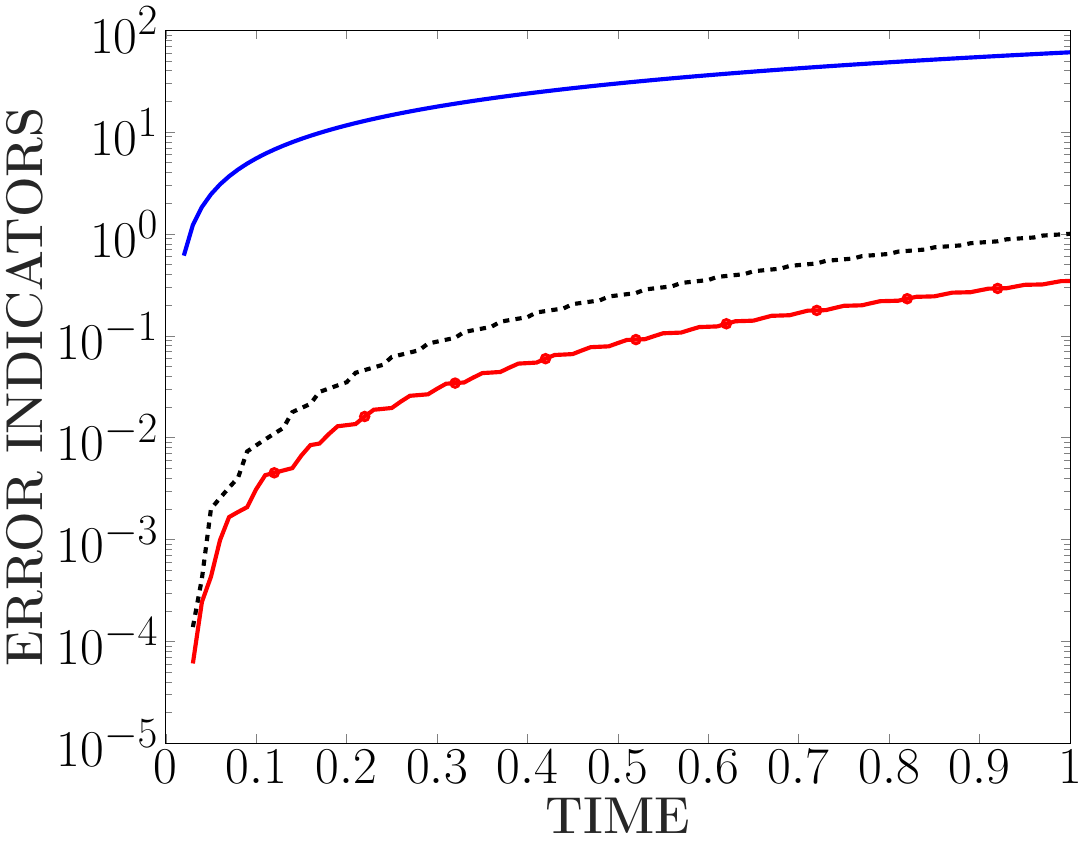}  
    \caption{Error indicator $\vartheta_1^n$ (solid line with circles) in \eqref{eq:time_error_indicator_II}, data approximation indicator $\delta^n$ (solid line) in \eqref{eq:data}, and LTS error indicator $\alpha^n$ (dotted line) in \eqref{eq:lts_error_indicator} vs. time.}
    \label{fig:r1-terms_accu_fixed}
\end{subfigure}
\end{figure}

In Fig. \ref{fig:solution_fixed}, we display the numerical solutions and the
underlying meshes for $h=0.3$ at times $t = 0.04$ and $t = 1$. The entire space-time time-evolving mesh
with $h^c = 0.3$ is shown in Fig. \ref{fig:mesh_fixed}. The refined
part moves to the right with the same unit speed as the Gaussian pulse. 
To verify convergence, we now consider a sequence of space-time meshes with
decreasing $h = h^c$ and time-step $\Delta t$, while keeping
all other parameters fixed.
Fig. \ref{fig:convergence_fixed} confirms that the LF-LTS-FEM method (\ref{eq:leapfrog}), including
local time-stepping, a nonzero source term, and a time-evolving mesh, indeed achieves the optimal
convergence rates $\mathcal{O}(h)$ and $\mathcal{O}(h^2)$ with
respect to the $\operatorname{H}^1(\Omega)$- and $\operatorname{L}^2(\Omega)$-norm, respectively.

\begin{figure}
\begin{subfigure}[t]{.475\textwidth}
  \centering
  \includegraphics[width=1\linewidth]{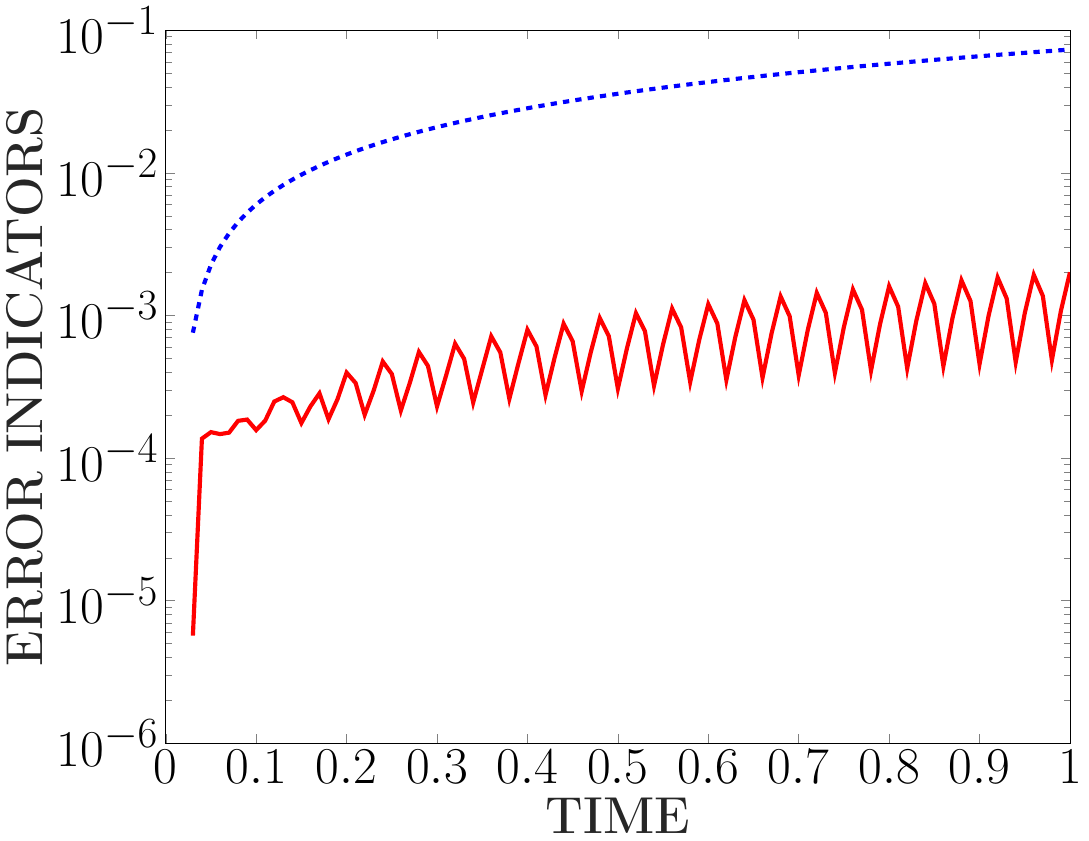}  
 \caption{Forced wave: Elliptic error indicator $\varepsilon_0^n$ in \eqref{eq:elliptic_error_indicator} (dotted line) and time error indicator $\vartheta_0^n$ in \eqref{eq:time_error_indicator_I} (solid line) vs. time without time accumulation.}
  \label{fig:r0-terms_fixed}
\end{subfigure}
\hfill
\begin{subfigure}[t]{.475\textwidth}
  \centering
  \includegraphics[width=1\linewidth]{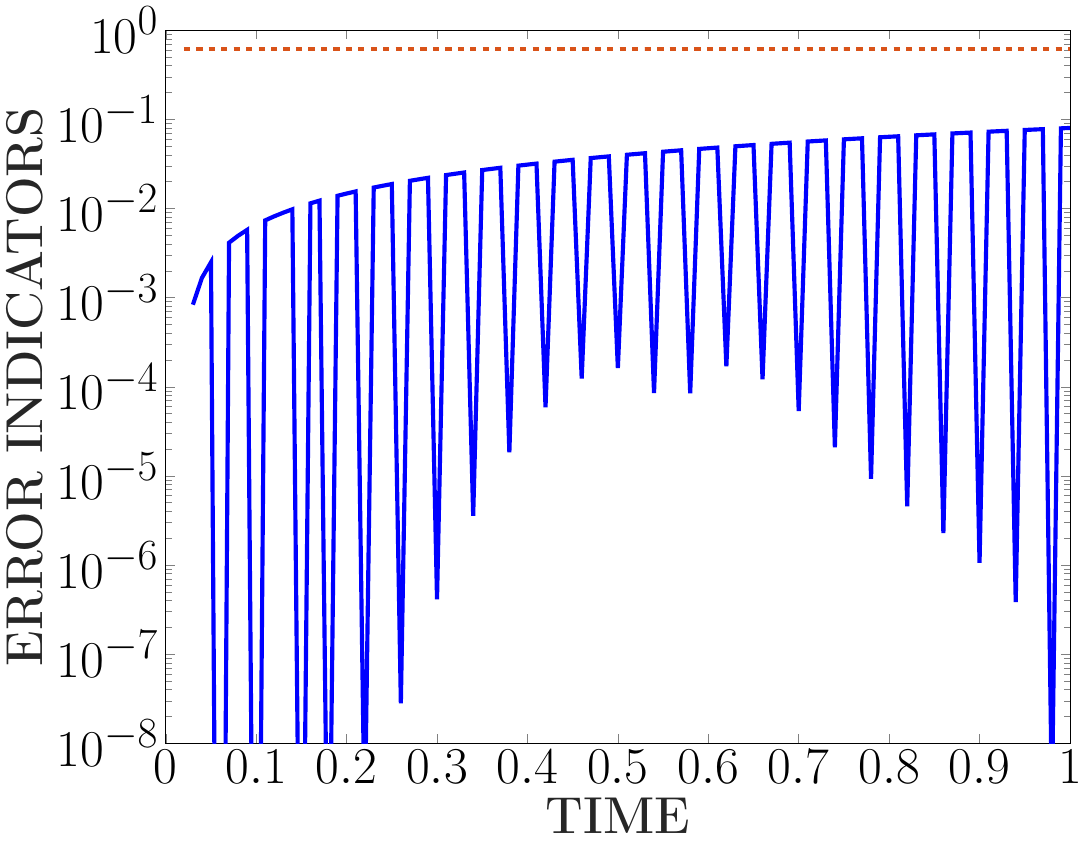}  
    \caption{Elliptic error indicator  $\varepsilon_1^n$ in \eqref{eq:elliptic_error_indicator} (solid line) and data approximation indicator $\delta^n$ (dotted line) in \eqref{eq:data} vs. time without time accumulation.}
    \label{fig:r1-terms_accu_fixed1}
\end{subfigure}
\newline
\begin{subfigure}[t]{.475\textwidth}
  \centering
  \includegraphics[width=1\linewidth]{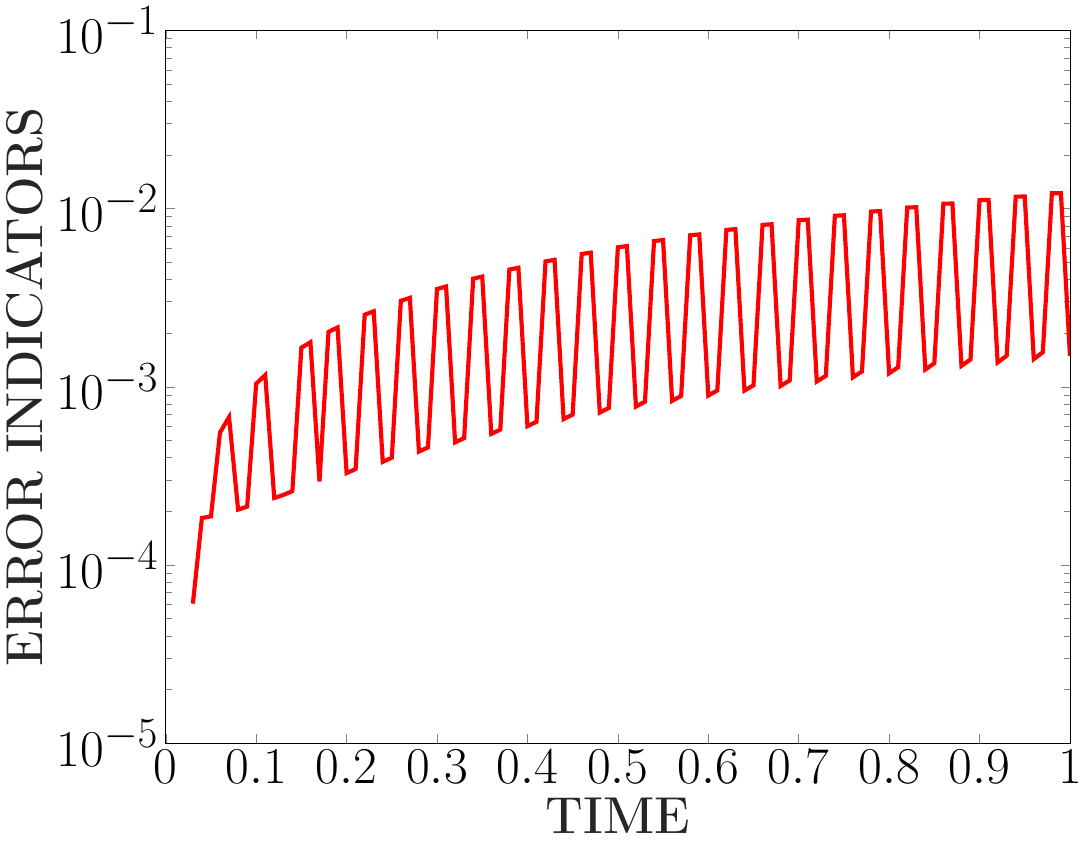}  
 \caption{Time error indicator $\vartheta_1^n$ in \eqref{eq:time_error_indicator_II} vs. time without time accumulation.}
  \label{fig:r1-terms_accu_fixed2}
\end{subfigure}
\hfill
\begin{subfigure}[t]{.475\textwidth}
  \centering
  \includegraphics[width=1\linewidth]{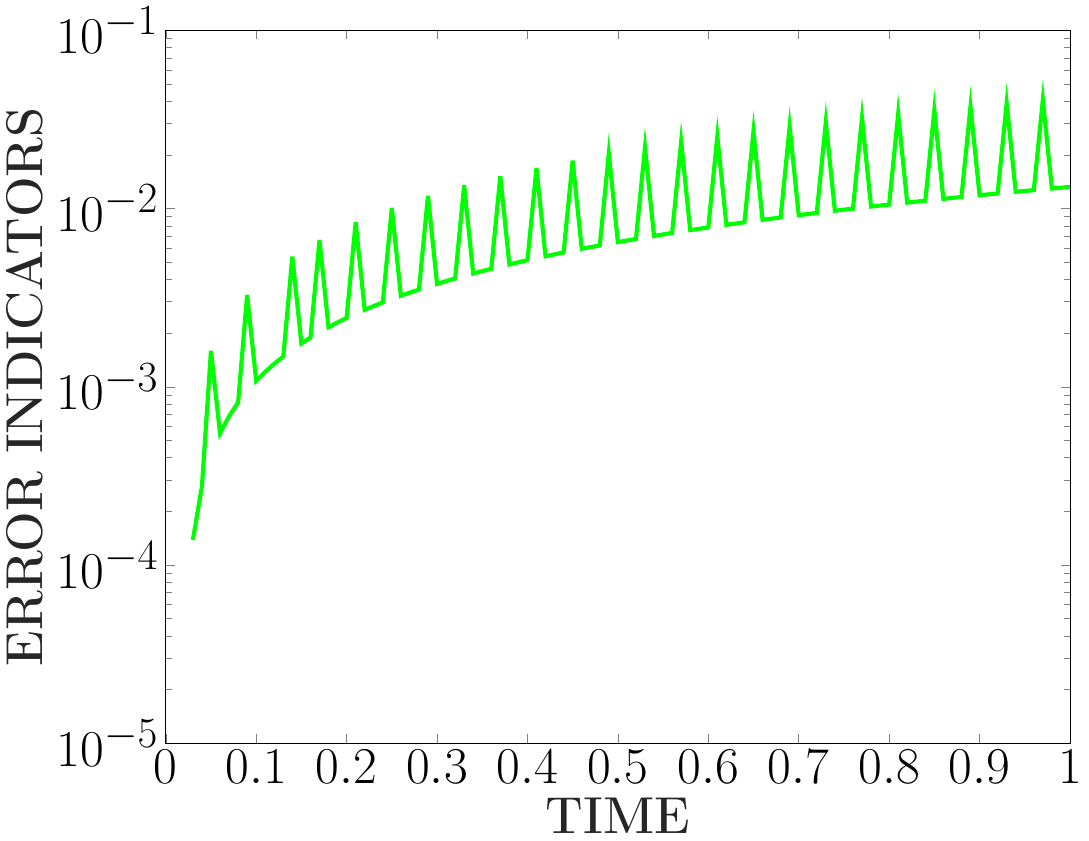}  
    \caption{LTS error indicator $\alpha^n$ \eqref{eq:lts_error_indicator} vs. time without time accumulation.}
    \label{fig:r1-terms_accu_fixed3}
\end{subfigure}
\end{figure}

Next, in Fig. \ref{fig:aposteriori_fixed}, we display the convergence rates of the full a posteriori
error estimates introduced in Section \ref{the:full-error analysis}. 
As expected, both converge as $\mathcal{O}(h)$ with the same rate as the numerical error with respect to the energy norm.
Further individual indicators from
Section \ref{the:full-error analysis} accumulated over time are shown
in Fig. \ref{fig:r0-terms_accu_fixed} and
Fig. \ref{fig:r1-terms_accu_fixed}.  The
behavior of the LTS error indicator $\alpha^n$ in
\eqref{eq:lts_error_indicator}, the time-error indicators
$\vartheta^n_0(t)$ and $\vartheta^n_1(t)$ and the data approximation indicator $\delta^n(t)$ together with the
elliptic error indicators $\varepsilon^n_0$ and $\varepsilon^n_1$ in
\eqref{eq:elliptic_error_indicator} are shown in Fig. \ref{fig:r0-terms_fixed}--\ref{fig:r1-terms_accu_fixed3} vs. time without accumulation. 
 As the space-time mesh is pre-defined and not adapted to the source $f(x,t)$, the data approximation indicator
 $\delta^n$ in Fig. \ref{fig:r1-terms_accu_fixed1} remains essentially constant over time.
The mesh-change indicators $\mu_0^n$ and $
\mu_1^n$ from Section \ref{eq:def:mesh-change-indicator} are not
displayed here, as mesh coarsening or refinement occurs only in regions
where the solution is nearly zero.

\subsection{Traveling wave}\label{sec:pulse}
Next, we apply the space-time adaptive algorithm from Section \ref{sect:algo}
to compute a rightward traveling wave; hence, the space-time mesh is no longer predefined but instead
automatically generated ''on the fly'' by the adaptive algorithm. 
The initial conditions $u_0, v_0$ define a right-moving Gaussian pulse with constant unit speed $c \equiv 1$ centered about $x = 1$ at $t = 0$:
\begin{equation}
  u(x,t) = e^{-4(x-1-t)^2}.
\end{equation}

\begin{figure}[t]
\begin{subfigure}[t]{.5\textwidth}
  \centering
  \includegraphics[width=1\linewidth,height=4.4cm]{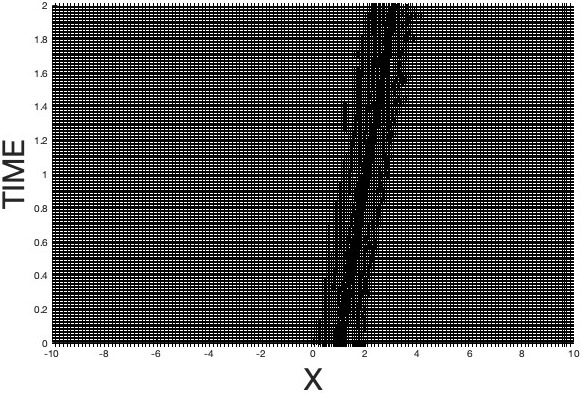}  
  \caption{Traveling wave: Time-adaptive mesh.}
  \label{fig:adaptive_mesh_pulse}
\end{subfigure}
\hfill
\begin{subfigure}[t]{.475\textwidth}
\centering
  \includegraphics[width=1\linewidth]{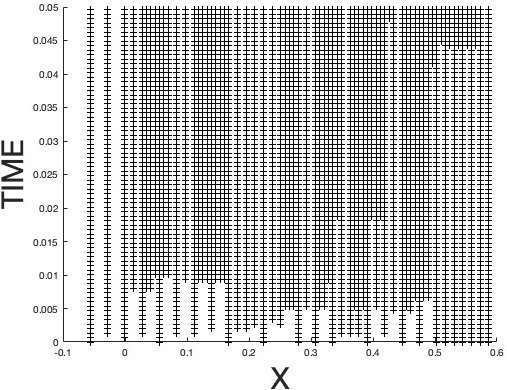} 
  \caption{Zoom of frame \ref{fig:adaptive_mesh_pulse} for $ t \in [0, 0.05] $ and $ x \in [-0.1, 0.6]$.}
  \label{fig:zoom_pulse}
\end{subfigure}
\newline
\begin{subfigure}[t]{.475\textwidth}
\centering
  \includegraphics[width=1\linewidth]{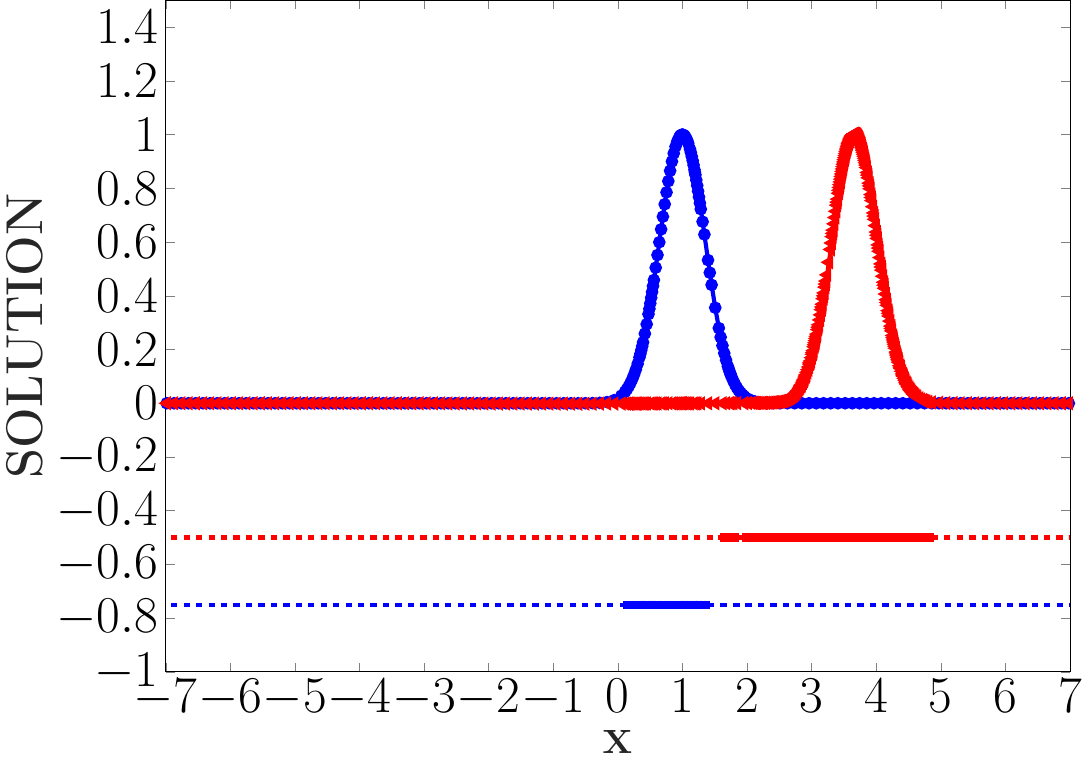}  
  \caption{Numerical solution and refined mesh at
time $t= 0$ (blue) and $t \approx 1.33$ (red).}
  \label{fig:solution_pulse}
\end{subfigure}
\hfill
\begin{subfigure}[t]{.475\textwidth}
\centering
  \includegraphics[width=1\linewidth]{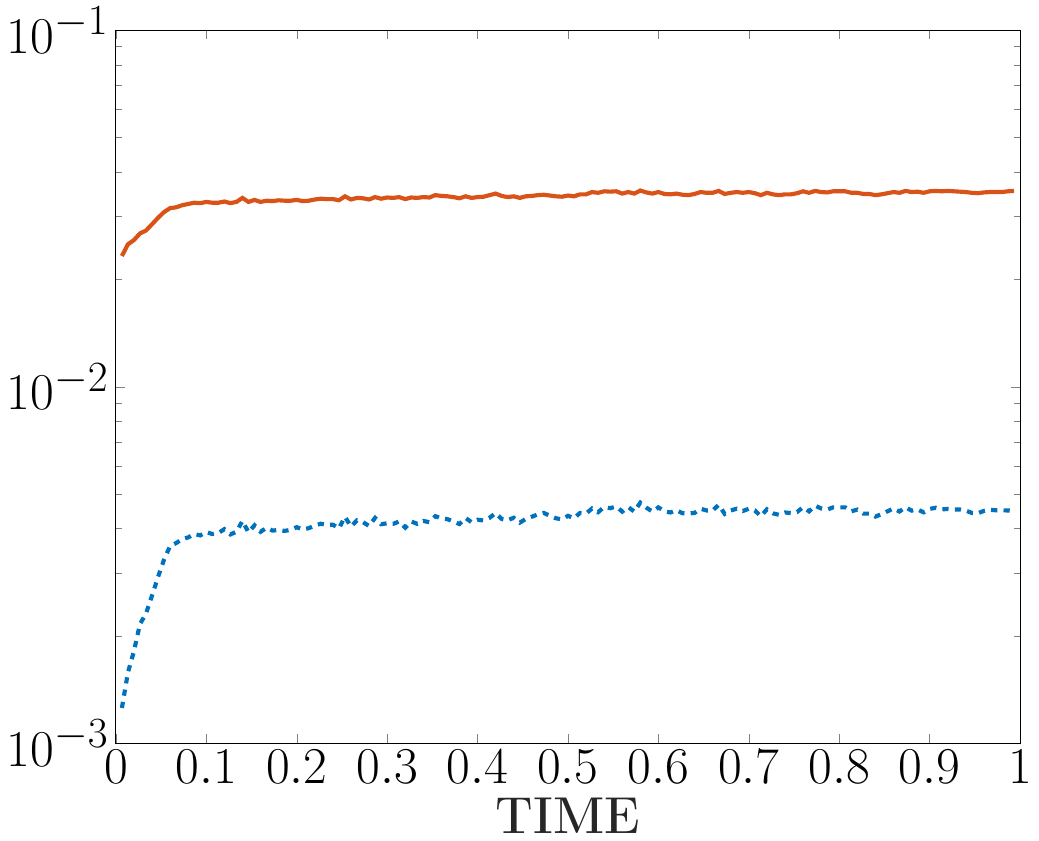}  
  \caption{Relative computational work (dotted line) and memory (solid line) w.r.t. uniform refinement.}
  \label{fig:memory_pulse}
\end{subfigure}
\end{figure}

  At any discrete time $t_n = n \; \Delta t$, the mesh $\mathcal{M}_n$ partitions the domain $\Omega$ into a coarse part, $\mathcal{M}_n^c$, of
  fixed mesh-size $h^c=h$, and a fine part, $\mathcal{M}_n^f$, of
  mesh-size $h^f = h^c / 2^k, \; k \leq 6$, where $k$ denotes the number of local refinements. Before time integration, the initially uniform coarse mesh $\mathcal{M}_0$ is adapted to the initial conditions $u_0, v_0$ by using \cref{alg:initial}. This yields the initial mesh $\mathcal{M}_1$ for the adaptive strategy. 
The global time-step $\Delta t$, proportional to $h^c$, corresponds to the CFL stability limit on a uniform mesh with mesh-size $h = h^c$.
To guarantee stability, the LF-LTS method locally adapts the 
time-step $\Delta \tau = \Delta t / p$ inside $\Omega^f_n$ proportionally to the local mesh size ratio $h^c / h^f$ -- see Section \ref{sec:lts}.
Hence the LF-LTS method takes $p = \left\lceil h_c / h_f \right\rceil$ local time-steps, at most $p = 2^6 = 64$, for each global time-step of size $\Delta t$ inside $\Omega^c_n$.

The entire space-time evolving mesh
with $h^c = 0.11$ is shown in Fig. \ref{fig:adaptive_mesh_pulse}.
We observe how the fine part of the mesh $\mathcal{M}_n^f$ automatically follows the peak of the Gaussian pulse as it propagates rightward across $\Omega$.  The mesh (and hence the associated FE space $\mathbbmss{V}_n$) changes whenever the maximum of the elliptic error indicators $\varepsilon^n_0$ and $ \varepsilon^n_1$ is larger than the given tolerance $\tol H$ divided by the total number of time-steps. Hence the fine mesh $\mathcal{M}^f_n$ moves to the right, with increasing $n$, at the same unit wave speed as the pulse, while two subsequent meshes $\mathbbmss{V}_n$ and $\mathbbmss{V}_{n+1}$ always remain compatible during any mesh change, see Section \ref{sec:compatible_mesh_change}. The remaining vertical lines inside the mesh are due to small dispersive waves in the discrete solution, which lead to localized refinement in regions of higher curvature. 

In Fig. \ref{fig:solution_pulse}, we display the numerical solutions and the underlying meshes for $h=0.5$ at times $t=0$ and $t \approx 1.33$.  Next, in Fig. \ref{fig:memory_pulse}, we compare both the amount of memory and the computational effort of the adaptive algorithm with a standard leapfrog FEM on a fixed mesh uniform both in space and time. The space-time adaptive LF-FEM algorithm requires at most $5\%$ of the number of dof's needed by a uniformly refined mesh with smallest mesh size $h = h^c/64$. To estimate the reduction in computational work, we compute the ratio of number of space-time dof's for
the adaptive LF-FEM algorithm vs. that using a standard scheme with uniform mesh size $h = h_c/64$ and time-step of $\Delta t/64$. The adaptive approach only requires
about $0.95\%$ of computational effort and thus achieves over a hundredfold reduction.

\begin{figure}
\begin{subfigure}[t]{.475\textwidth}
  \centering
  \includegraphics[width=1\linewidth]{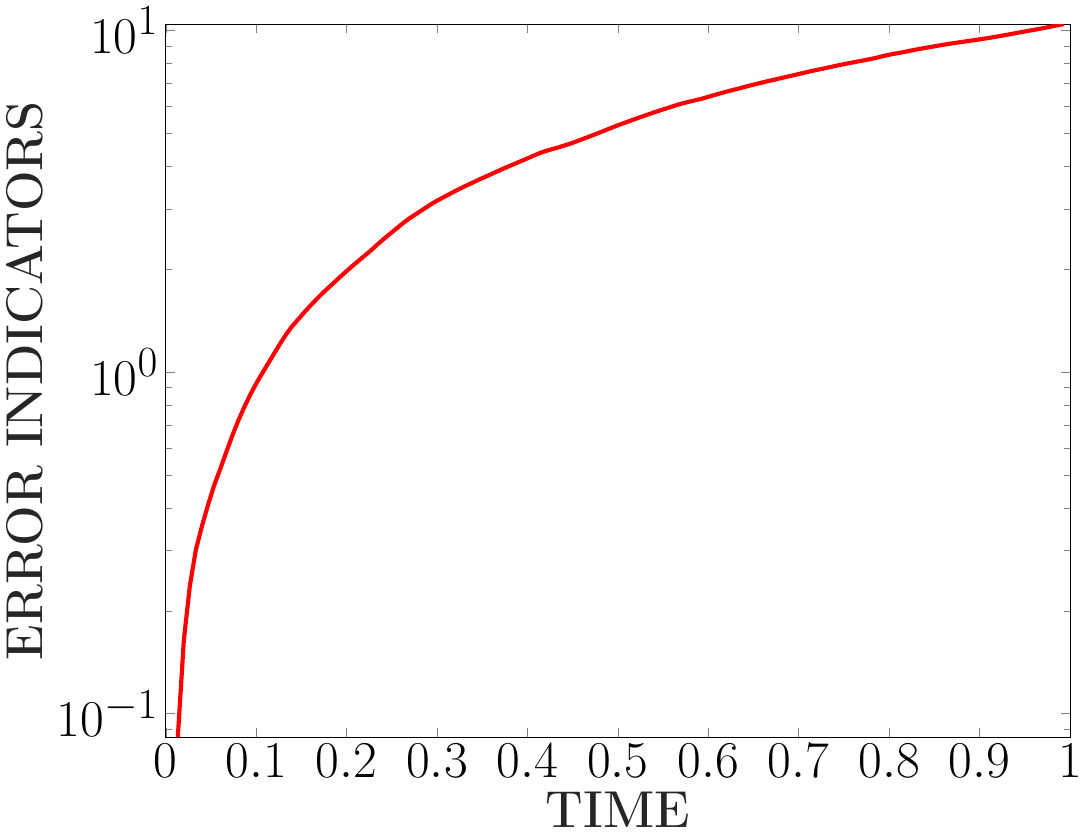}  
 \caption{Traveling wave: Time evolution of the error indicator $\vartheta_0^n$ in \eqref{eq:lts_error_indicator}.}
  \label{fig:r0-terms_accu_pulse}
\end{subfigure}
\hfill
\begin{subfigure}[t]{.475\textwidth}
  \centering
  \includegraphics[width=1\linewidth]{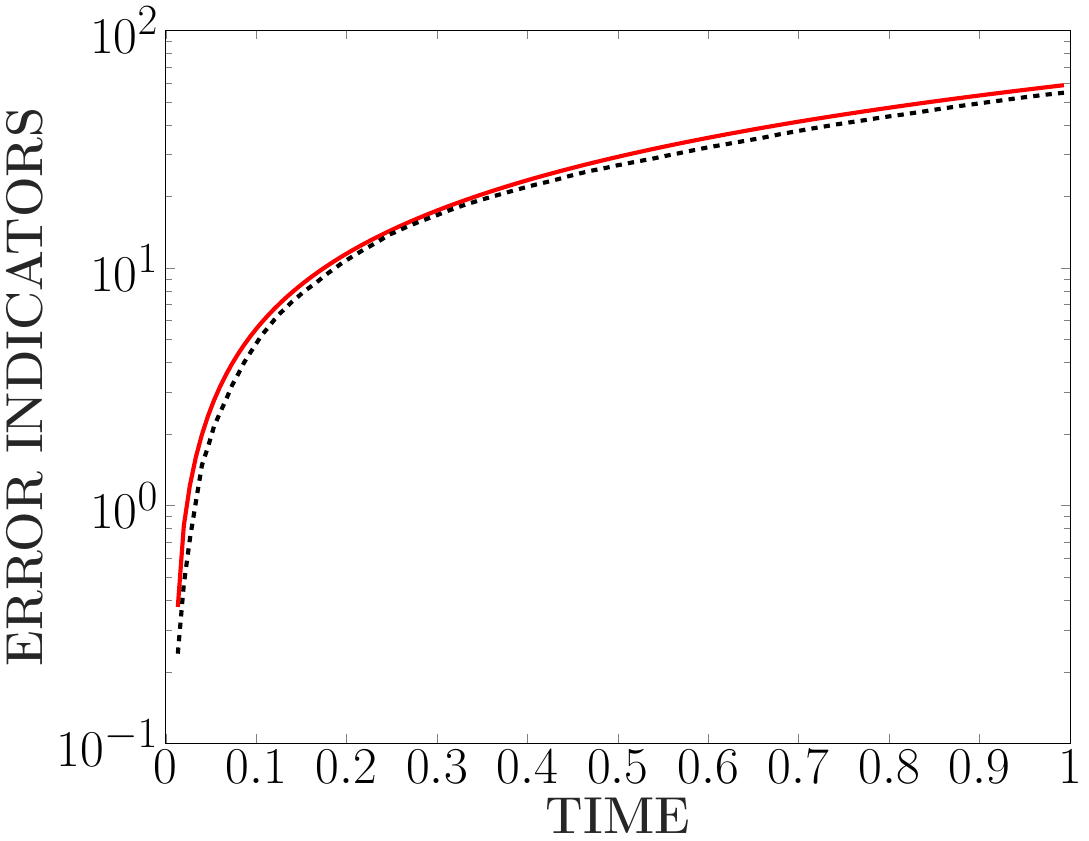}  
    \caption{Time evolution of the time error indicator $\vartheta_1^n$ (solid line) in \eqref{eq:time_error_indicator_II} and the LTS error indicator $\alpha^n$ (dashed line) in \eqref{eq:lts_error_indicator}.}
    \label{fig:r1-terms_accu_pulse}
\end{subfigure}
\newline
\begin{subfigure}[t]{.475\textwidth}
  \centering
  \includegraphics[width=1\linewidth]{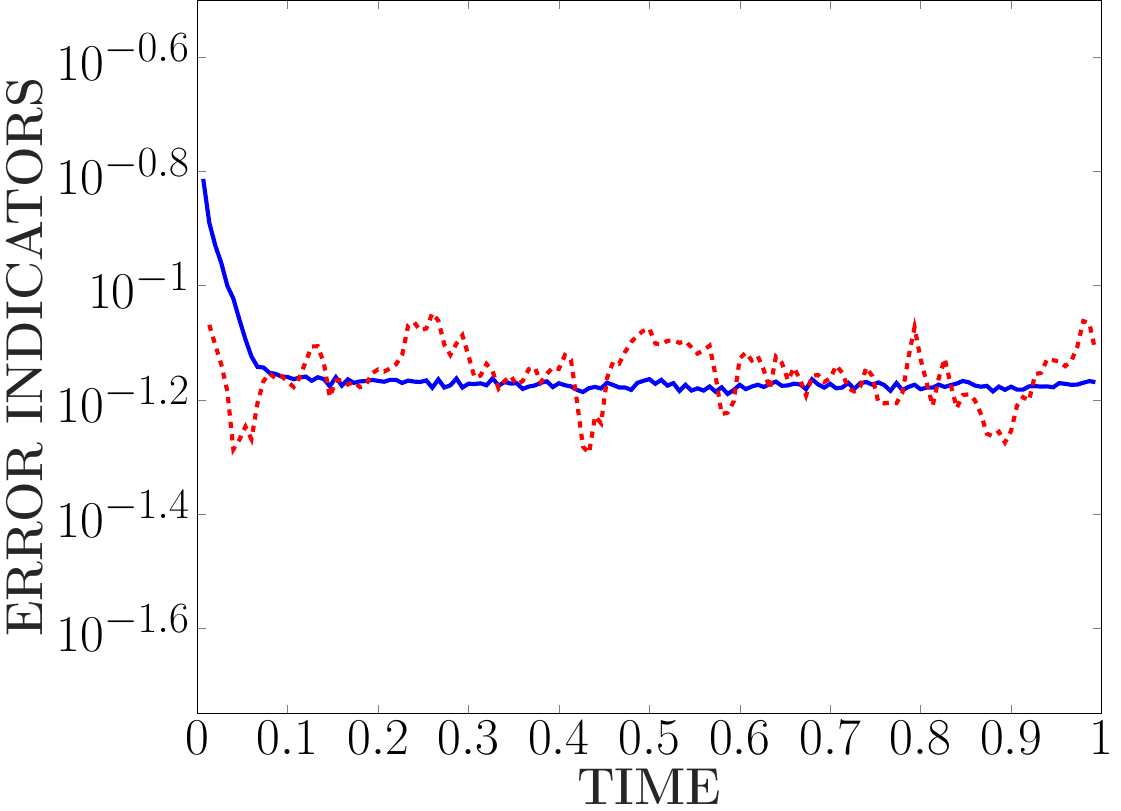}  
 \caption{Elliptic error indicator $\varepsilon_0^n$ in \eqref{eq:elliptic_error_indicator} (solid line) and time error indicator $\vartheta_0^n$ in \eqref{eq:time_error_indicator_I} (dashed line) vs. time without time accumulation.}
  \label{fig:r0-terms_pulse}
\end{subfigure}
\hfill
\begin{subfigure}[t]{.475\textwidth}
  \centering
  \includegraphics[width=1\linewidth]{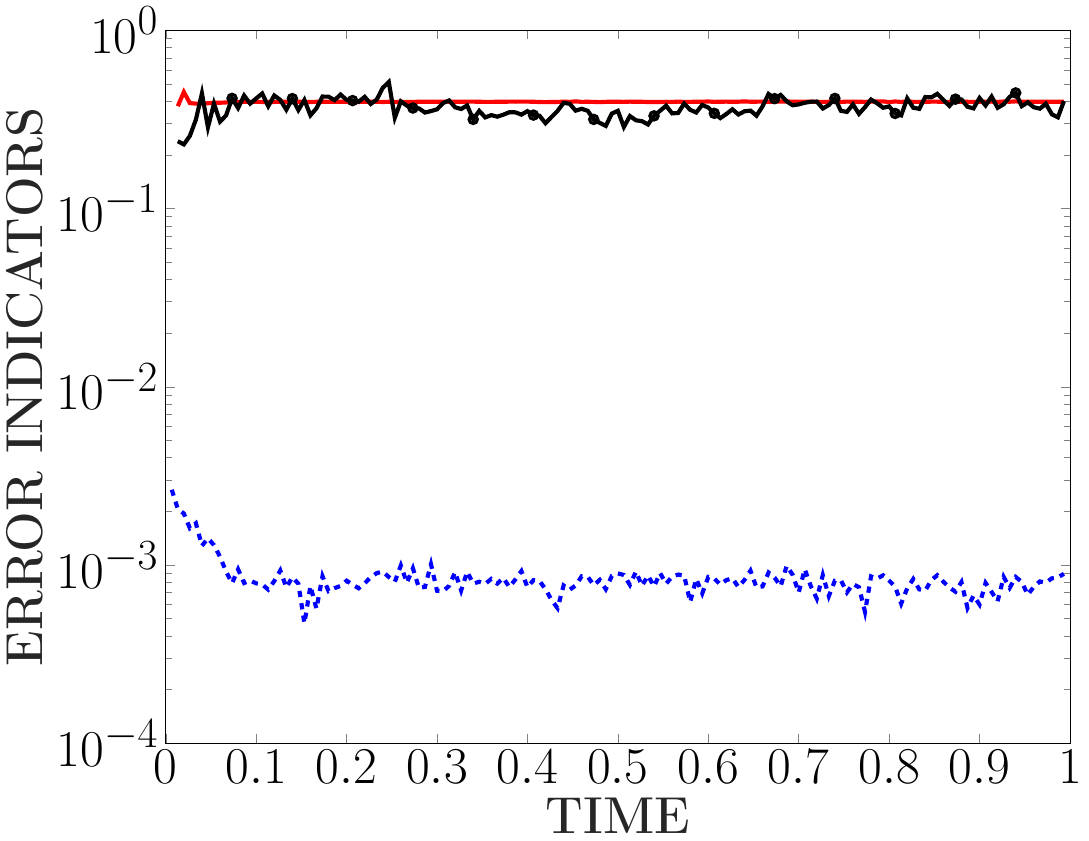}  
    \caption{Elliptic error indicator  $\varepsilon_1^n$ in \eqref{eq:elliptic_error_indicator} (dashed line), time error indicator $\vartheta_1^n$ (solid line) in \eqref{eq:time_error_indicator_II}, and LTS error indicator $\alpha^n$ (solid line with circles) in \eqref{eq:lts_error_indicator} vs. time without time accumulation.}
    \label{fig:r1-terms_pulse}
\end{subfigure}
\end{figure}

In Fig. \ref{fig:r0-terms_accu_pulse} and Fig. \ref{fig:r1-terms_accu_pulse}, we display the time evolution of  various a posteriori error indicators from
Section \ref{the:full-error analysis} accumulated over time. The
behavior of the LTS error indicator $\alpha^n$ in
\eqref{eq:lts_error_indicator} and time-error indicators
$\vartheta_0^n(t)$ and $\vartheta_1^n(t)$ in \eqref{eq:time_error_indicator_I} together with the
elliptic error indicators $\varepsilon_0^n$ and $\varepsilon_1^n$ in \eqref{eq:elliptic_error_indicator} are shown in Fig. \ref{fig:r0-terms_pulse} and Fig. \ref{fig:r1-terms_pulse} vs. time without accumulation. The mesh-change indicators $\mu_0^n$ and $\mu_1^n$ \eqref{eq:def:mesh-change-indicator} are not displayed here, as mesh coarsening/refinement occurs only in regions where the solution is nearly zero, hence they remain vanishingly small. Since the source $f$ is identically zero, the data approximation indicator $\delta^n(t)$ also remains identically zero in this example.

\subsection{Splitting wave}
Next, we consider
a Gaussian pulse which splits into two symmetric waves traveling in opposite directions with constant speed $c \equiv 1$.
The initial conditions $u_0, v_0$ are set such that exact solution, centered about $x = 1$ at $t = 0$, is
\begin{equation}
  u(x,t) = \frac{1}{2}\left[e^{-4(x-1-t)^2} + e^{-4(x-1+t)^2} \right].
\end{equation}

\begin{figure}
\begin{subfigure}[t]{.5\textwidth}
  \centering
  \includegraphics[width=1\linewidth,height=4.2cm]{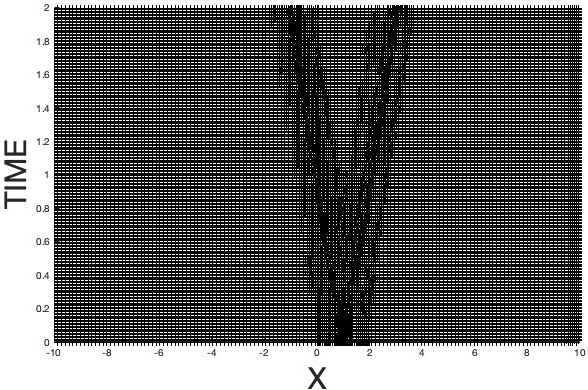}  
  \caption{Splitting wave: Space-time adaptive mesh.}
  \label{fig:adaptive_mesh_splitting}
\end{subfigure}
\hfill
\begin{subfigure}[t]{.45\textwidth}
\centering
  \includegraphics[width=1\linewidth]{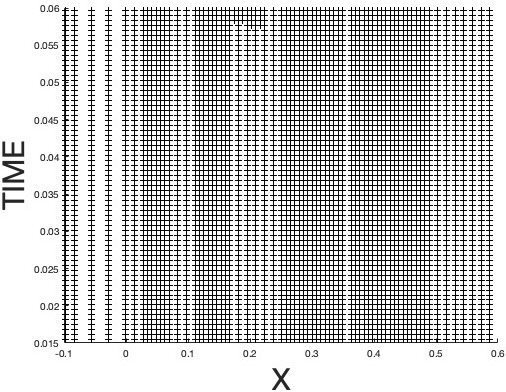}  
  \caption{Zoom of frame \ref{fig:adaptive_mesh_splitting} for $t\in [0.015, 0.06]$ and $x\in [-0.1, 0.6]$.}
  \label{fig:zoom_splitting}
\end{subfigure}
\newline
\begin{subfigure}[t]{.475\textwidth}
\centering
  \includegraphics[width=1\linewidth]{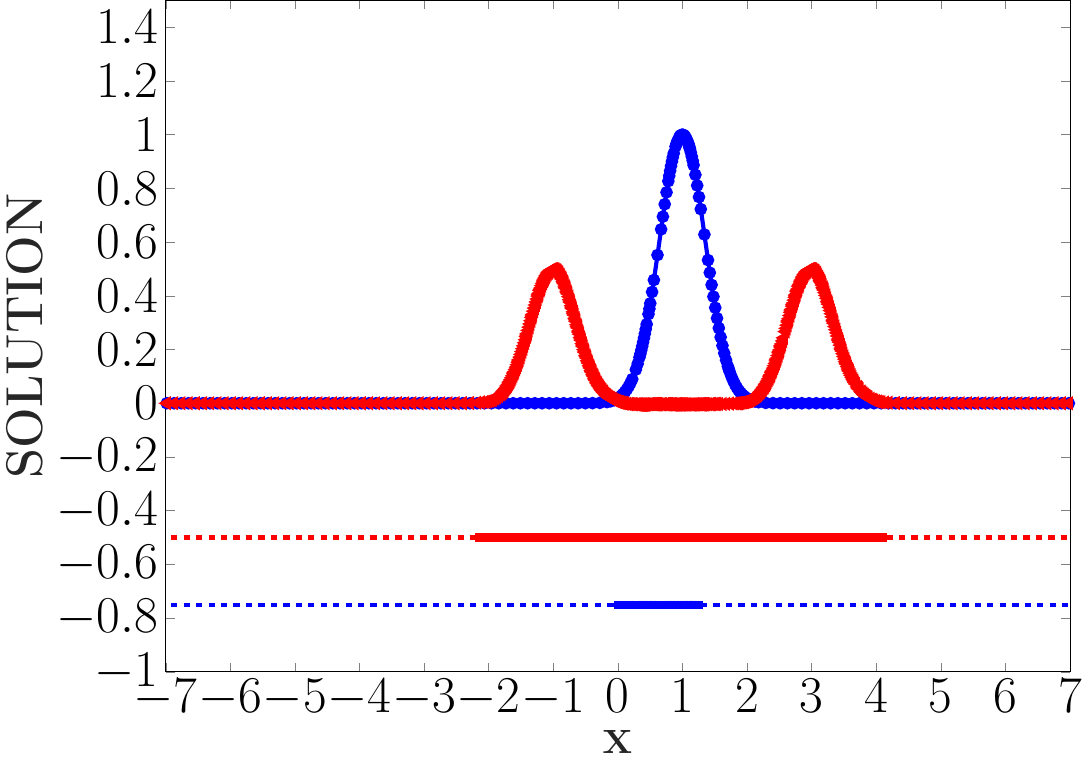}  
  \caption{Numerical solution and refined mesh at
time $t= 0$ (blue) and $t \approx 1.33$ (red).}
  \label{fig:solution_splitting}
\end{subfigure}
\hfill
\begin{subfigure}[t]{.475\textwidth}
\centering
  \includegraphics[width=1\linewidth]{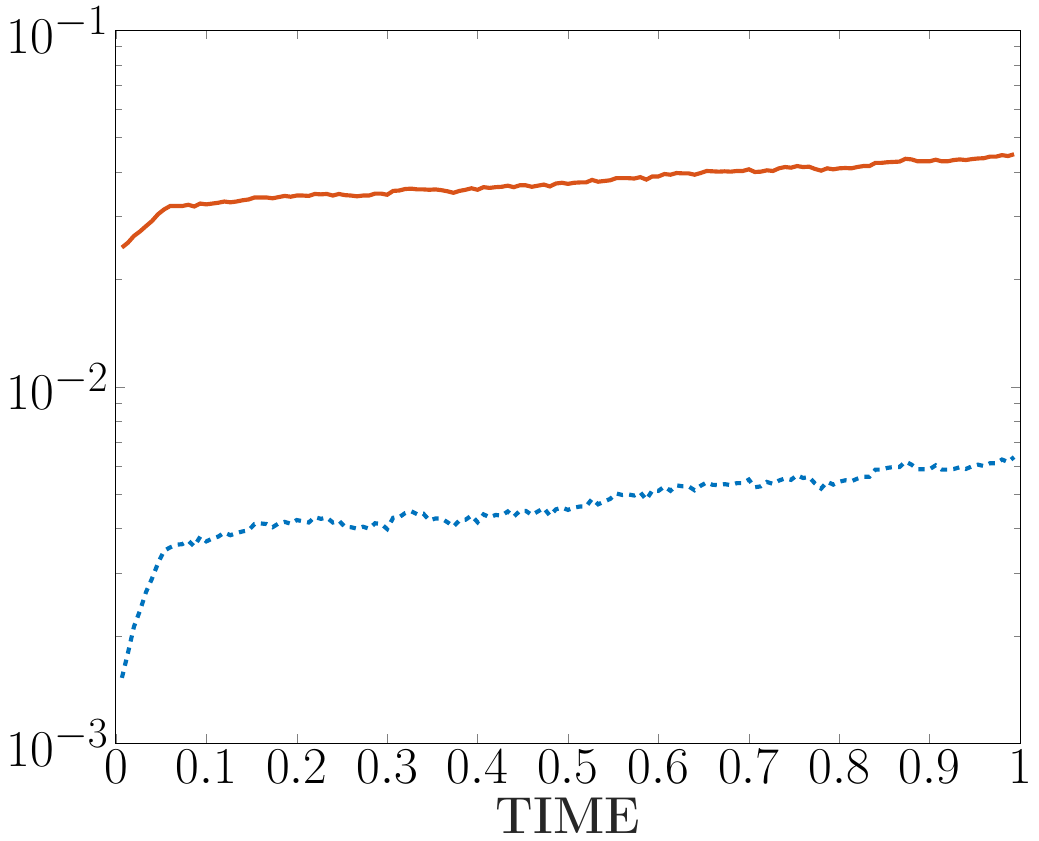}  
  \caption{Relative computational work (dashed line) and memory (solid line) w.r.t. uniform refinement.}
  \label{fig:memory_splitting}
\end{subfigure}
\end{figure}

Again we apply the space-time 
  adaptive algorithm from Section \ref{sect:algo} as described in Section \ref{sec:pulse}. We observe how the fine part of the mesh $\mathcal{M}_n^f$ automatically adapts and follows the two separating peaks of the Gaussian pulse. In particular, the initially refined single interval automatically splits into two separate refined regions, each associated with its own local time-step. That topological change in the refined part of the mesh requires no particular attention or reordering of the unknowns, as the time integration is fully explicit.
  The mesh (and hence the associated FE space $\mathbbmss{V}_n$) again changes whenever the maximum of the elliptic error indicators $\varepsilon^n_0$ and 
  $\varepsilon^n_1$ is larger than the given tolerance $\tol H$ divided by the total number of time-steps. During any mesh change,  two subsequent meshes $\mathbbmss{V}_n$ and $\mathbbmss{V}_{n+1}$ always remain compatible, as shown in Fig. \ref{fig:adaptive_mesh_splitting} for $h^c = 0.11$.
  
In Fig. \ref{fig:solution_splitting}, we display the numerical solutions and the underlying meshes for $h=0.5$ at the initial time $t=0$ and at
time $t \approx 1.33$. Next, in Fig. \ref{fig:memory_splitting}, we compare both the amount of memory and the computational effort of the adaptive algorithm with a standard leapfrog FEM on a fixed mesh uniform both in space and time. The space-time adaptive LF-FEM algorithm requires at most $6\%$ of the number of dof's needed by a uniformly refined mesh with mesh-size $h = h^c/64$. 
 To estimate the reduction in computational work, we compute at every time-step the ratio of the number of space-time dof's for
the adaptive LF-FEM algorithm vs. that using a standard scheme with uniform mesh size $h = h_c/64$ and time-step of $\Delta t/64$. The adaptive approach only requires
about $1.2\%$ of computational effort and thus again
achieves close to a hundredfold reduction.

\begin{figure}
\begin{subfigure}[t]{.475\textwidth}
  \centering
  \includegraphics[width=1\linewidth]{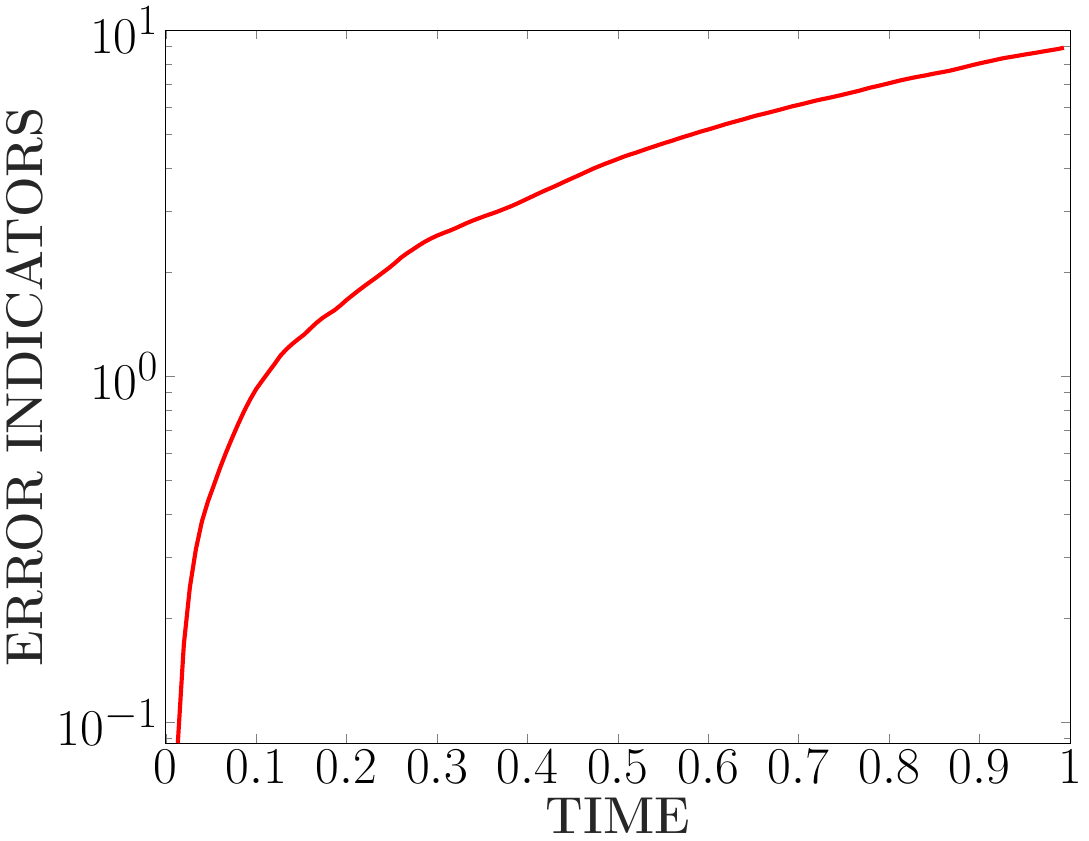}  
 \caption{Splitting wave: Time evolution of the error indicator $\vartheta_0^n$ in \eqref{eq:lts_error_indicator}.}
  \label{fig:r0-terms_accu}
\end{subfigure}
\hfill
\begin{subfigure}[t]{.475\textwidth}
  \centering
  \includegraphics[width=1\linewidth]{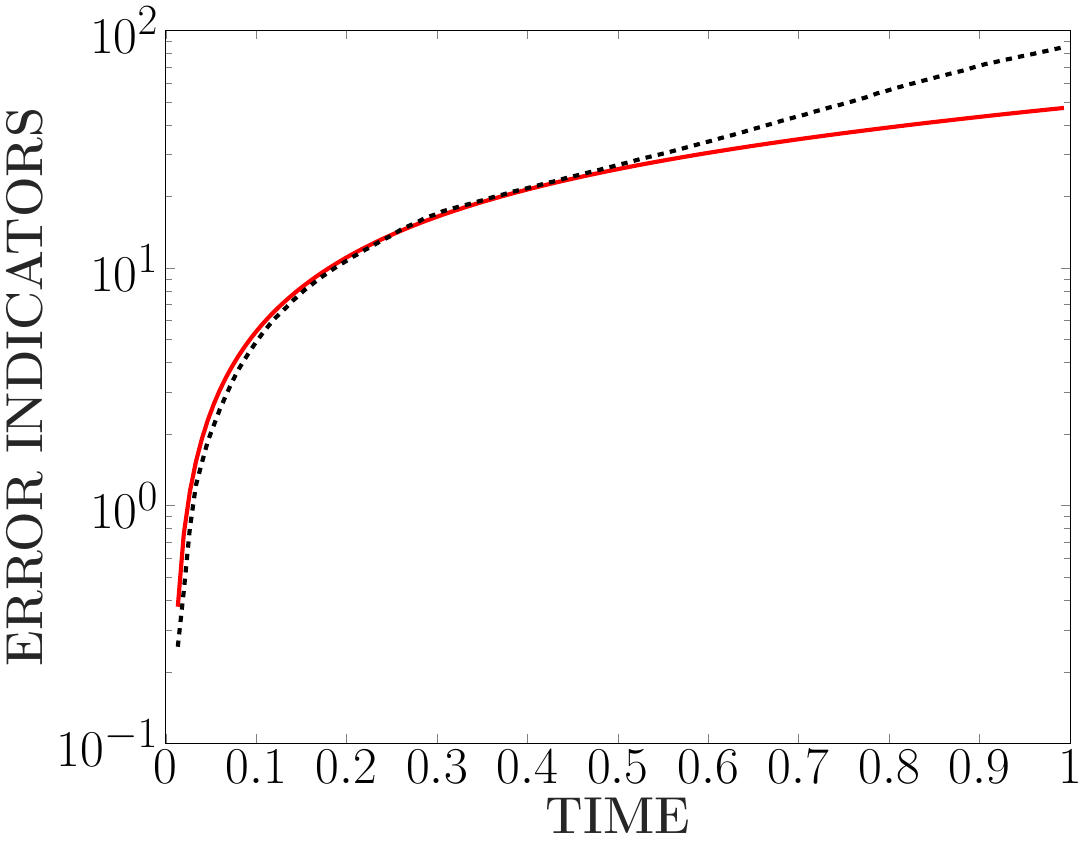}  
    \caption{Time evolution of the time error indicator $\vartheta_1^n$ (solid line) in \eqref{eq:time_error_indicator_II} and the LTS error indicator $\alpha^n$ (dashed line) in \eqref{eq:lts_error_indicator}.}
    \label{fig:r1-terms_accu}
\end{subfigure}
\newline
\begin{subfigure}[t]{.475\textwidth}
  \centering
  \includegraphics[width=1\linewidth]{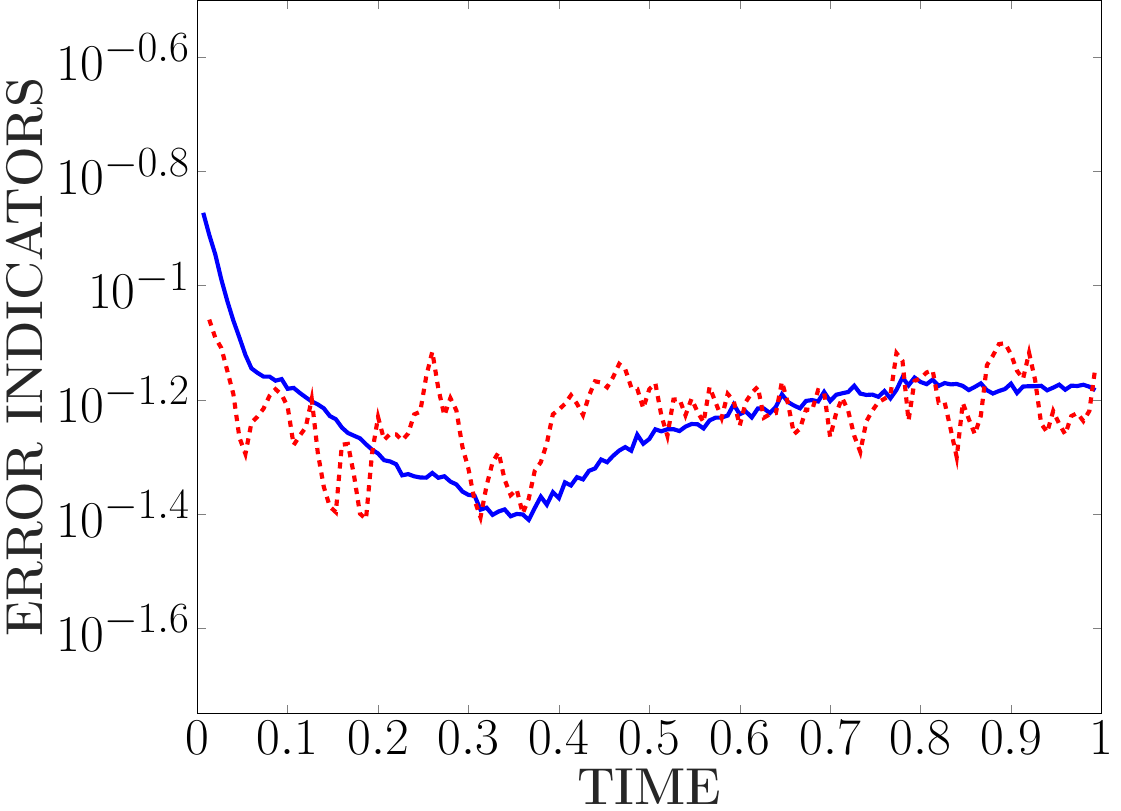}  
 \caption{Elliptic error indicator $\varepsilon_0^n$ in \eqref{eq:elliptic_error_indicator} (solid line) and time error indicator $\vartheta_0^n$ in \eqref{eq:time_error_indicator_I} (dasehd line) vs. time without time accumulation.}
  \label{fig:r0-terms}
\end{subfigure}
\hfill
\begin{subfigure}[t]{.475\textwidth}
  \centering
  \includegraphics[width=1\linewidth]{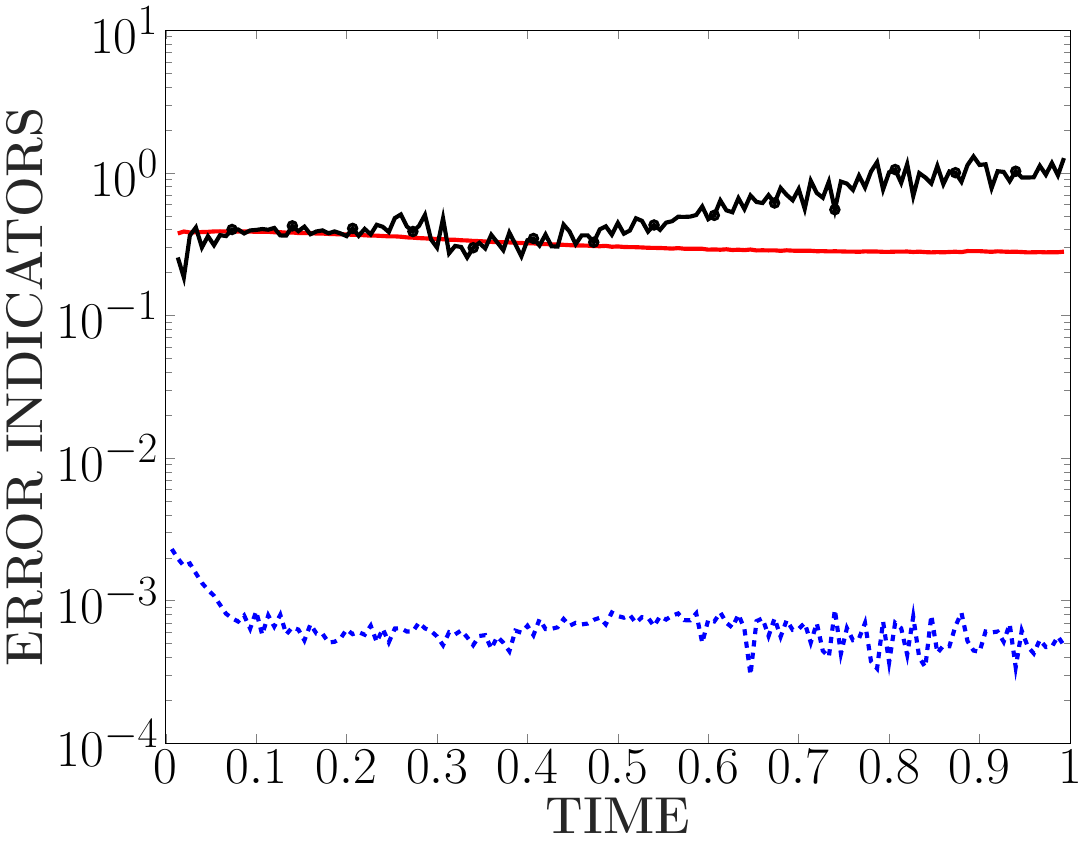}  
    \caption{Elliptic error indicator  $\varepsilon_1^n$ in \eqref{eq:elliptic_error_indicator} (dashed line), time error indicator $\vartheta_1^n$ (solid line) in \eqref{eq:time_error_indicator_II}, and LTS error indicator $\alpha^n$ (solid line with circles) in \eqref{eq:lts_error_indicator} vs. time without time accumulation.}
    \label{fig:r1-terms}
\end{subfigure}
\end{figure}

In Fig. \ref{fig:r0-terms_accu} and Fig. \ref{fig:r1-terms_accu}, we display various a posteriori error indicators from Section \ref{the:full-error analysis} accumulated over time. The behavior of the LTS error indicator $\alpha^n$ in \eqref{eq:lts_error_indicator} and time-error indicators $\vartheta_0^n(t)$ and $\vartheta_1^n(t)$ together with the
elliptic error indicators $\varepsilon_0^n$ and $\varepsilon_1^n$ in \eqref{eq:elliptic_error_indicator} in Fig. \ref{fig:r0-terms} and Fig. \ref{fig:r1-terms} decay with time. Again, the mesh-change indicators $\mu_0^n$ and $\mu_1^n$ \eqref{eq:def:mesh-change-indicator} are not displayed here, as they remain vanishingly small. Since the source $f$ is identically zero, the data approximation indicator $\delta^n(t)$ also vanishes here.

\subsection{L-shaped domain}
Finally, we consider a Gaussian pulse which propagates with constant wave speed $c \equiv 1$ across an L-shaped domain.
The initial data is chosen such that at $t = 0$ the solution is a Gaussian
centered at $(x,y) = (0.4,\,0.6)$:
\begin{equation}
  u_0(x,y) = e^{-600\left((x-0.4)^2 + (y-0.6)^2\right)},
\end{equation}
while the initial velocity $v_0 \equiv 0$.
At all boundaries we impose homogeneous Dirichlet boundary conditions
and choose a uniform initial mesh $\mathcal M_0$ with $h=0.08$.
Again we apply the space–time adaptive algorithm from Section~\ref{sect:algo}, following the set-up described in Section~\ref{sec:pulse}. The refinement and coarsening are performed using the iFEM package \cite{Chen:2008ifem}, combined with our stabilized LF–LTS method for time integration with a fixed stabilization parameter $\nu=0.01 $.

\begin{figure}
\begin{subfigure}[t]{.475\textwidth}
  \centering
  \includegraphics[width=1\linewidth]{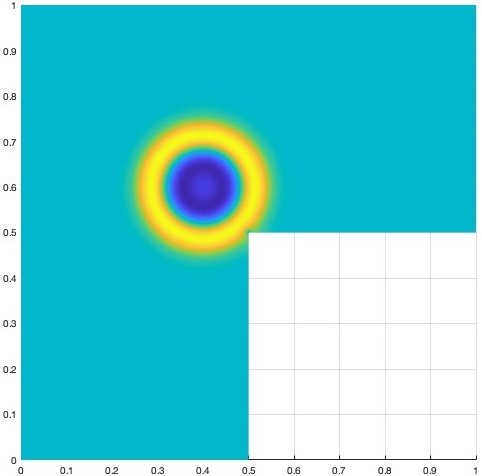}  
 \caption{L-shaped domain: Numerical solution at
time $t= 0.1$.}
  \label{fig:2DWave_01}
\end{subfigure}
\hfill
\begin{subfigure}[t]{.475\textwidth}
  \centering
  \includegraphics[width=1.2 \linewidth,height=5.7cm]{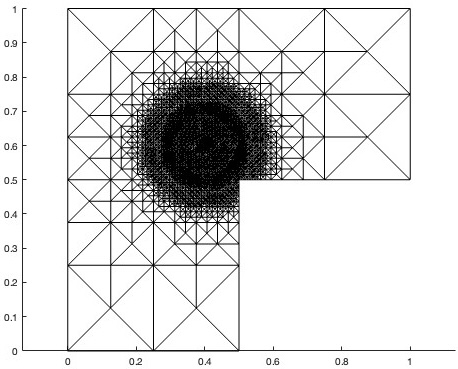} 
    \caption{Space-time adaptive mesh at
time $t= 0.1$.}
    \label{fig:2DWaveMesh_01}
\end{subfigure}
\newline
\begin{subfigure}[t]{.475\textwidth}
  \centering
  \includegraphics[width=1\linewidth]{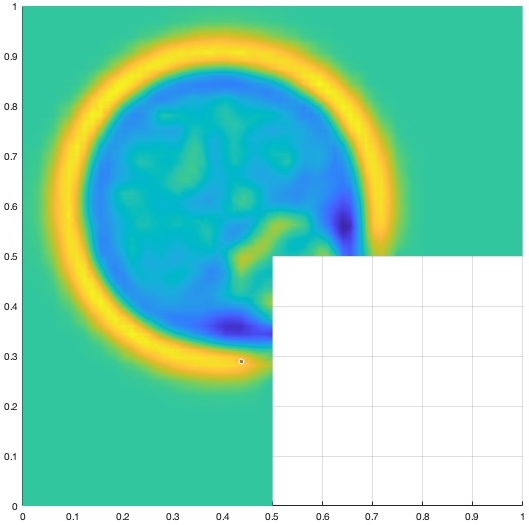}  
 \caption{Numerical solution at
time $t= 0.3$.}
  \label{fig:2DWave_03}
\end{subfigure}
\hfill
\begin{subfigure}[t]{.475\textwidth}
  \centering
  \includegraphics[width=1.2\linewidth,height=5.7cm]{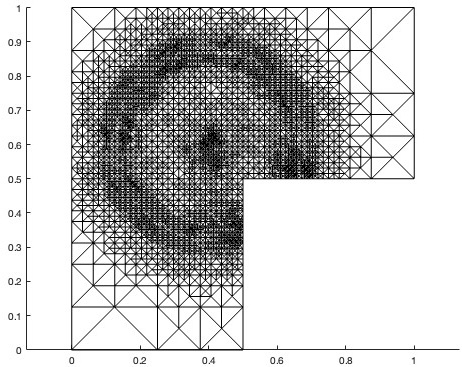} 
    \caption{Space-time adaptive mesh at
time $t= 0.3$.}
    \label{fig:2DWaveMesh_03}
\end{subfigure}
\end{figure}

\begin{figure}
\begin{subfigure}[t]{.475\textwidth}
  \centering
  \includegraphics[width=1\linewidth,height=4.7cm]{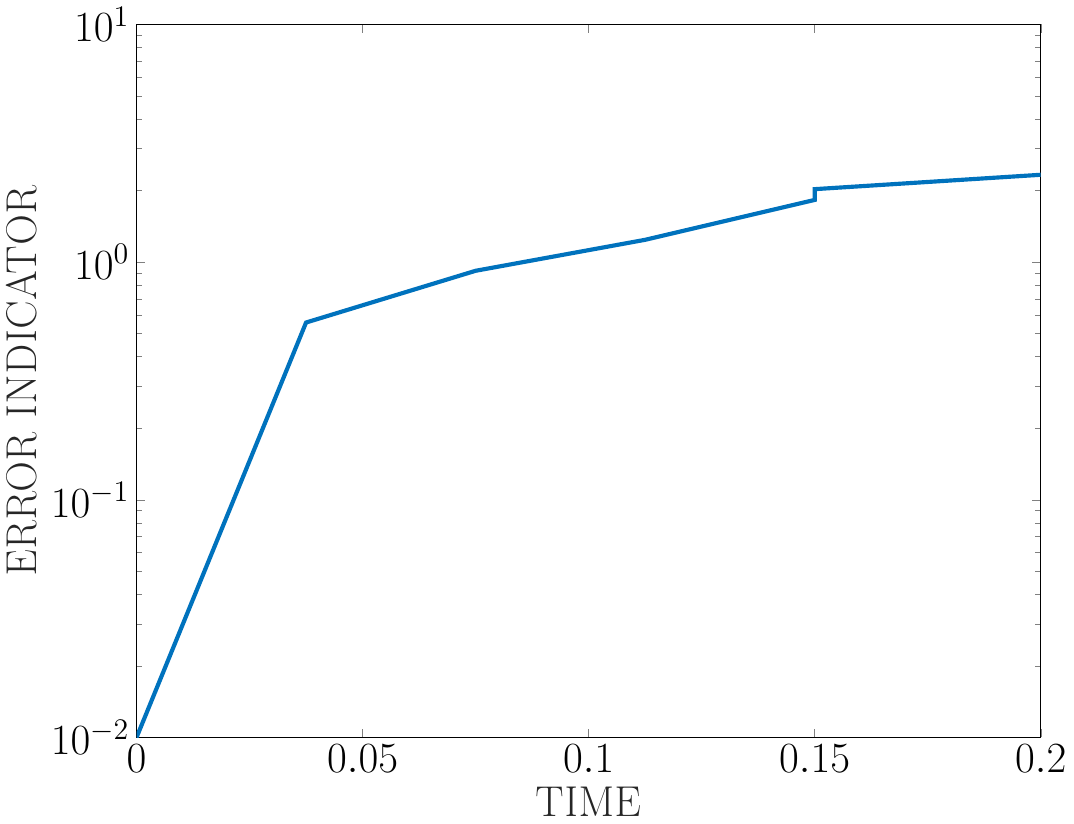}  
 \caption{L-shaped domain: Elliptic error indicator $\varepsilon^n_0$ in \eqref{eq:elliptic_error_indicator} vs.\ time without accumulation.}
  \label{fig:estimate2D}
\end{subfigure}
\hfill
\begin{subfigure}[t]{.475\textwidth}
  \centering
  \includegraphics[width=1\linewidth]{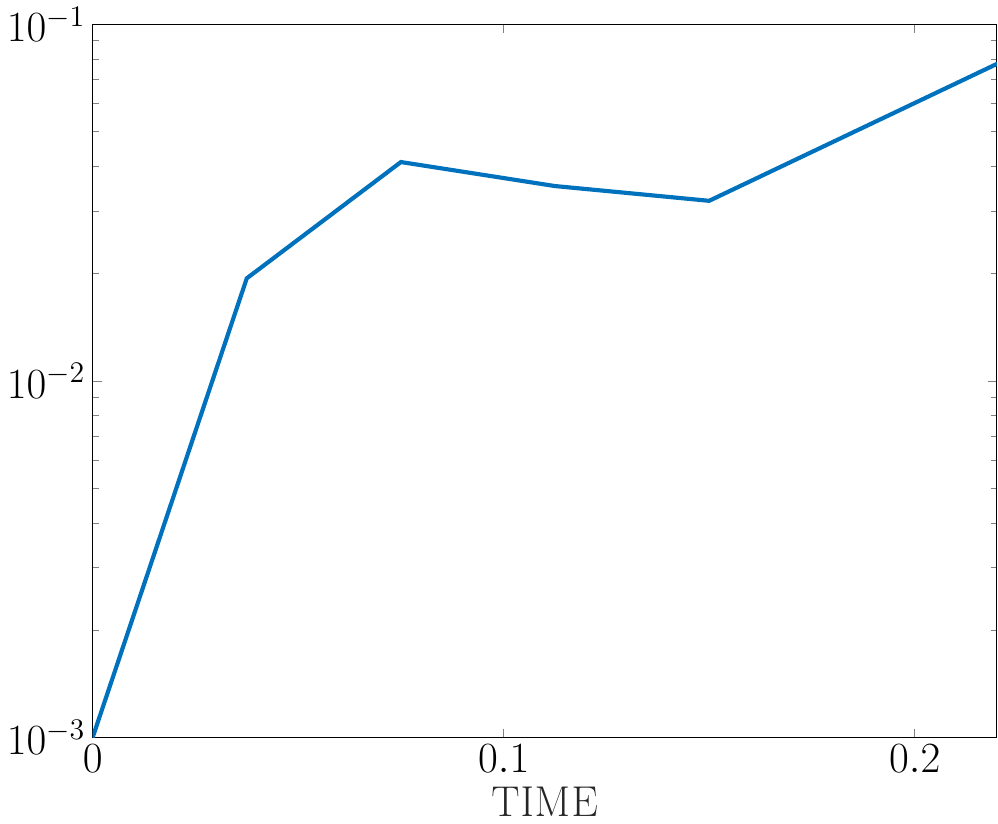} 
    \caption{Relative memory w.r.t. uniform refinement.}
    \label{fig:memory2D}
\end{subfigure}
\end{figure}

As shown in Figures~\ref{fig:2DWave_01} and \ref{fig:2DWave_03}, the refined region of the mesh $ \mathcal{M}_n^f$ automatically adapts and tracks the propagating wave front.
The corresponding adaptive meshes are shown in Figures~\ref{fig:2DWaveMesh_01} and \ref{fig:2DWaveMesh_03}.
In Fig. \ref{fig:estimate2D}, we display the elliptic error indicator $\varepsilon_0^n$ \eqref{eq:elliptic_error_indicator} from Section \ref{the:full-error analysis} without accumulation over time. Finally, in Fig.\ \ref{fig:memory2D}, 
we show the relative memory requirement of the adaptive LF-FEM by computing 
the ratio of the number of dof's (FE nodes in the mesh) for
the adaptive LF-FEM algorithm vs. that using a uniform mesh.
\section{Concluding remarks}
We have proposed a space-time adaptive strategy for the solution of the time-dependent wave equation which combines three essential ingredients for robustness and efficiency. First, the rigorous a posteriori estimates \cite{GroteLakkisSantos:24:PosterioriErrorEstimatesWave} from Section 3.1 provide guaranteed error bounds that include all sources of error. As the 
error indicators are fully computable and local, they permit to decide "on the fly" whether to accept, or recompute on an adapted mesh, the numerical solution at the current time-step. 
Second, by restricting mesh change to compatible meshes, as 
described in Section 3.2, we ensure that local mesh refinement 
never leads to a loss in accuracy while allowing coarsening only
where appropriate to minimize the resulting inherent information loss
from interpolation or projection. Third, we apply leapfrog based
explicit local time-stepping  \cite{DiazGrote:09:article:Energy-Conserving, 
GroteMitkova:10:article:Explicit-Local,
  GroteMehlinSauter:18:article:Convergence,
  GroteMichelSauter:21:article:Stabilized}
inside any locally refined region of the mesh, each associated with its own local time-step. Thus we overcome the bottleneck caused
by local mesh refinement due to any overly stringent CFL stability condition 
while retaining the simplicity, efficiency and inherent parallelism of explicit time integration. 

Although our adaptive algorithm relies on a fixed underlying 
coarse mesh and a constant global time-step, it nonetheless automatically adapts both the mesh and the time-step to capture locally the waves' smaller scale features as they propagate across the computational domain. 
While the increased flexibility of a varying global time-step might be all too tempting, one must keep in mind that the standard leapfrog 
method with varying time-step is prone to instability \cite{Skeel:93:article}.

Our numerical experiments confirm the expected order of convergence and illustrate the usefulness of our adaptive strategy. 
Although it is only studied in one and two space dimensions,
our space-time adaptive strategy, including the a posteriori error bounds and local time-stepping
approach, are dimension independent and immediately apply three space dimensions, too. 
Yet even in a single space dimension,
the computational effort and memory requirement are up to two 
orders of magnitude smaller than a standard approach on a uniform mesh.
Clearly, that reduction in computational cost and memory requirement will be even more pronounced in higher dimensions, although
operations related to mesh change then also become more involved.
The overall gain in efficiency compared to a static, uniform space-time mesh will generally depend on the problem at hand and in particular on the size of
the locally refined region relative to the entire space-time mesh. 
Nonetheless
the reduction in the number of degrees of freedom will generally enable simulations of improved accuracy, which otherwise would be prohibitive merely due to sheer problem size.

\section*{Acknowledgements} This work was supported by the Swiss National Science Foundation under grant SNF 200020-188583.

\printbibliography
\end{document}